# Globally Optimal Energy-Efficient Power Control and Receiver Design in Wireless Networks

Alessio Zappone, *Senior Member, IEEE*, Emil Björnson, *Member, IEEE*, Luca Sanguinetti, *Senior Member, IEEE*, and Eduard Jorswieck, *Senior Member, IEEE*

***Abstract*—The characterization of the global maximum of energy efficiency (EE) problems in wireless networks is a challenging problem due to the non-convex nature of investigated problems in interference channels. The aim of this work is to develop a new and general framework to achieve globally optimal solutions. First, the hidden monotonic structure of the most common EE maximization problems is exploited jointly with fractional programming theory to obtain globally optimal solutions with exponential complexity in the number of network links. To overcome this issue, we also propose a framework to compute suboptimal power control strategies characterized by affordable complexity. This is achieved by merging fractional programming and sequential optimization. The proposed monotonic framework is used to shed light on the ultimate performance of wireless networks in terms of EE and also to benchmark the performance of the lower-complexity framework based on sequential programming. Numerical evidence is provided to show that the sequential fractional programming framework achieves global optimality in several practical communication scenarios.**

## I. Introduction

The percentage of the global footprint (in terms of $CO_2$-equivalent emissions) due to the information and communications technology (ICT) was recently estimated to be 5% [1]. Although this is a small percentage, it is rapidly increasing, and the situation will escalate in the near future with the advent of 5G networks. Credited sources foresee the number of connected devices to reach 50 billions by 2020 [2] and that the data traffic will increase by 1000× over the next 10 years [3]. If no countermeasures are taken, the energy demand to operate and provide such massive data rates to so many devices will become unmanageable, and the resulting greenhouse gas emissions and electromagnetic pollution will exceed safety thresholds. While restricting the global ICT usage is unrealistic, a promising answer to this issue lies in optimizing the energy efficiency (EE) of ICT systems; that is, in maximizing the amount of information reliably transmitted per Joule of consumed energy. Moreover, the EE is of paramount importance for operators (e.g., to save on electricity bills) and for end-users (e.g., to prolong the lifetime of batteries). This has created a great interest in the design of new network architectures, as well as new beamforming and power control strategies taking into account the cost of energy so as to push the EE in communication systems towards new limits.

A. Zappone is with the University of Cassino and Southern Lazio, Cassino, Italy (alessio.zappone@unicas.it). E. Björnson is with the Department of Electrical Engineering (ISY), Linköping University, Linköping, Sweden (emil.bjornson@liu.se). L. Sanguinetti is with the University of Pisa, Dipartimento di Ingegneria dell'Informazione, Pisa, Italy (luca.sanguinetti@unipi.it) and also with the Large Systems and Networks Group (LANEAS), CentraleSupélec, Gif-sur-Yvette, France. E. Jorswieck is with TU Dresden, Faculty of Electrical and Computer Engineering, Communications Laboratory, Dresden, Germany (eduard.jorswieck@tu-dresden.de).

The work of A. Zappone was supported by the German Research Foundation (DFG), under grant CEMRIN - ZA 747/1-3.

E. Björnson was supported by the Swedish Foundation for Strategic Research (SFF).

L. Sanguinetti was supported by the ERC Starting Grant 305123 MORE and by the research project 5GIOTTO funded by the University of Pisa.

The work of E. Jorswieck was supported by the German Research Foundation (DFG), within the Collaborative Research Center 912 HAEC.

A preliminary version of this paper will be presented at ICASSP 2016, Shanghai, China, 20-25 March 2016.

### A. State-of-the-art

The EE of a wireless communication link is commonly defined as a benefit-cost ratio, where the data rate is compared with the associated energy consumption:

$$\text{EE [bit/Joule]} = \frac{\text{Data rate [bit/s]}}{\text{Power Consumption [W]}}.$$

Given the fractional nature of the EE, the main mathematical tool for the optimization of EE-related metrics is fractional programming [4]—a branch of optimization theory that provides algorithms with polynomial complexity to globally maximize fractional functions with a concave numerator and a convex denominator [5]. However, even this powerful tool fails when interference-limited networks must be optimized. This is due to the fact that the presence of multi-user interference makes the numerator of the EE a non-concave function of the transmit power. A common way to circumvent this problem is to only consider suboptimal orthogonal or semi-orthogonal transmission schemes as well as interference cancellation techniques, to fall back into the noise-limited case. Contributions in this direction are given in [6]–[8]. In [6], [7], multi-carrier networks are considered and the global energy efficiency (GEE) of the system (defined as the ratio between the achievable sum rate and the total power consumption) is optimized using orthogonal or semi-orthogonal subcarrier allocation schemes. In [8], the authors consider a multi-antenna system and aim at maximizing the GEE when non-linear interference cancellation techniques are used. In [9], the presence of a large number of base station antennas is also exploited to average out multi-user interference, whereas in [10] energy efficiency optimization for systems using wireless information and power transfer is performed, always considering the regime of large deployed antennas. In [11] zero-forcing is considered in the context of wireless power transfer. However, interference suppression techniques as well as orthogonal transmissions inevitably require a large amount of radio resources to be implemented, and thus are not practical in large networks. Besides, unavoidable channel estimation errors also prevent from achieving perfect interference suppression.

Alternative available approaches typically try to handle the interference by embedding alternating optimization techniques into fractional programming methods. Examples in this context can be found in [12] wherein the minimum of the individual EEs is maximized, in [13] where the maximization of GEE is considered, and in [14] which considers the maximization of the sum of the GEE of the different base stations of the network. While these approaches can be applied to interference-limited networks (and are also able to tackle beamforming problems in the case of [13], [14]), they can not ensure the convergence to resource allocation points enjoying strong optimality properties. Moreover, they are tailored to specific problem formulations, considering specific instances of communication systems.

A first attempt to provide a unified framework to tackle the optimization of general EE metrics in a centralized way is based on the integration of traditional fractional programming methods with the tool of sequential optimization [15]–[19]. The basic idea of sequential optimization is to tackle a difficult problem by solving a sequence of easier approximate problems, which can be solved by standard methods. Provided that suitable approximations can be found, sequential optimization is able to obtain a solution, which fulfills first-order optimality conditions for the original problem, while at the same time requiring only the solution of convex problems. Sequential optimization was first used for rate and signal-to-interference-plus-noise ratio (SINR) optimization in [15], [16], and it was more recently successfully integrated into fractional programming to tackle EE maximization problems. Contributions in this sense are [17]–[20], which consider both multi-antenna and multi-carrier systems, also addressing two-hop communications and full-duplex systems. However, sequential optimization can not guarantee to achieve global optimality, and indeed the above works do not provide any insight into the gap between sequential optimization algorithms and the globally optimal solution. In general, no previous work provides efficient and provably convergent algorithms to obtain the global maximum of the EE in interference-limited networks.

A possible answer to close this gap is represented by monotonic optimization. This optimization framework can solve optimization problems where the utility and constraints are monotonically increasing/decreasing functions of the variables, even if no convexity assumption is fulfilled. The polyblock algorithm [21] and the branch-reduce-and-bound (BRB) algorithm [22] are two key algorithms for globally solving these problems. These algorithms have recently been used to solve different non-convex problems in communication and networking systems; for example, joint power control and scheduling [23] and sum rate maximization with multi-antenna transmitters [24]. The weighted sum-rate maximization problem is studied by using monotonic optimization and the rate profile technique in [25] for Gaussian interference channels. In [26], monotonic optimization is used to characterize the Pareto boundary of rate optimization problems in multi-cell networks. A good overview of monotonic optimization techniques in communication and networking is provided by [27] and [28]. However, to-date monotonic optimization has never been used for EE optimization, the main difficulty being that the EE is not a monotone function of the transmit powers, which is the case with more traditional performance metrics like rate or SINR.

### B. *Major contributions*

In contrast to previous related literature, this work aims at shedding some light on the ultimate EE performance of interference-limited networks. This overall goal is accomplished by making the following major contributions:

- Two novel optimization frameworks for interference-limited networks are developed: the former obtains global optimality, whereas the latter fulfills first-order optimality conditions, with a reduced complexity. Both frameworks encompass a wide class of energy-efficient maximization problems. More specifically:
  1) Both frameworks are able to work with multiple energy-efficient performance metrics, such as the GEE, the minimum of the individual energy efficiencies, and the product of the individual energy efficiencies. Additionally, the globally optimal framework is also able to handle the sum of the individual energy efficiencies of network nodes.
  2) Both frameworks are able to handle multiple types of constraints, such as maximum power constraints, minimum rate constraints, and interference temperature constraints. Both frameworks are able to perform both transmit power allocation, and joint transmit power and receiver design optimization. Moreover, the extension to beamforming optimization problems is also possible and briefly discussed.
  3) Both frameworks encompass multiple expressions of the links' SINRs, which generalize the typical SINR expression commonly encountered in previous related literature. Considering multiple SINR models makes the proposed frameworks broad enough to be applied to many relevant instances of contemporary and future wireless communication networks, such as massive MIMO networks, heterogeneous networks, LTE networks, device-to-device communications, full-duplex networks.
- Global optimality is obtained by merging fractional programming and monotonic optimization, developing a new *monotonic fractional programming* framework. Due to the use of monotonic optimization, the complexity of the proposed framework turns out to be exponential in the number of links, but still lower than standard global optimization algorithms. Moreover, convergence to the global optimum is theoretically ensured, whereas this is not always true for general global optimization methods.
- The monotonic fractional programming framework is used to get insights on the ultimate EE performance of wireless

networks. The energy-efficient Pareto boundary of the network is characterized. Moreover, monotonic fractional programming provides an effective benchmark for any low complexity, yet suboptimal, resource allocation algorithm.

• Fractional programming is used together with sequential optimization to develop a novel *sequential fractional programming* framework able to obtain candidate solutions fulfilling the Karush-Kuhn-Tucker (KKT) optimality conditions of EE maximization problems. The framework has limited complexity, requiring only to solve convex optimization problems. The effectiveness of such low-complexity solutions is validated by a numerical analysis which employs the monotonic fractional solution as a benchmark. The results indicate that the proposed sequential fractional solution achieves near-optimal performance.

### C. Outline and notation

The remainder of this paper is organized as follows. Section II defines the signal model and formulates the EE maximization problems. Section III introduces some useful mathematical definitions and results on fractional programming, monotonic optimization, and sequential programming. Sections IV and V develop the monotonic and sequential fractional programming frameworks, respectively. Numerical results for two notable case-studies of communication systems are illustrated in Section VI, whereas concluding remarks are made in Section VII.

The following notation is used throughout the paper. Scalars are denoted by lower case letters whereas boldface lower case letters are used for vectors. The superscripts $^T$ and $^H$ denote transpose and conjugate transpose, respectively. $\mathbf{1}_N$ and $\mathbf{0}_N$ are the $N$-dimensional all-one and all-zero vectors, respectively. $\mathbb{R}$ denotes the real number space and $\mathbb{C}$ is the complex number space. $\mathbb{R}^N$ stands for the $N \times 1$ real vector space and $\mathbb{R}_+^N$ denotes its non-negative orthant. $\mathbb{R}_{++}$ denotes the set of strictly positive real numbers. $\nabla_{\mathbf{x}} f$ denotes the gradient vector of function $f(\mathbf{x})$ with respect to $\mathbf{x}$ and $|\mathbf{A}|$ stands for the determinant of a matrix $\mathbf{A}$. For $\mathbf{x} \in \mathbb{R}^N$ and $\mathbf{y} \in \mathbb{R}^N$, we use $\mathbf{x} \succeq \mathbf{y}$ to indicate that $\mathbf{x}$ is greater than or equal to $\mathbf{y}$ in a component-wise manner.

## II. Signal model and Problem formulation

Consider a wireless network wherein $K$ mutually interfering links are active over a communication bandwidth $B$ [in Hz]. Each link includes a single-antenna transmitter node and a receiver node (possibly equipped with multiple antennas). Call $p_k$ the transmit power level [in W] of link $k$ (from the transmitter to its intended receiver) and assume that $0 \leq p_k \leq P_{\max,k}$ where $P_{\max,k}$ is the maximal transmit power.

Denote by $\gamma_k(\mathbf{p})$ the SINR of link $k$ as a function of $\mathbf{p} = [p_1, \ldots, p_K] \in \mathbb{R}_+^K$. At this stage, no particular expression for the function $\gamma_k(\mathbf{p})$ will be specified while only the following general assumption is made:

**Assumption 1.** *The function $\gamma_k(\mathbf{p}) : \mathbb{R}_+^K \to \mathbb{R}_+$ is, for all $k$, such that the achievable rate $R_k(\mathbf{p})$ of link $k$ can be expressed as the difference of two non-negative functions:*

$$R_k(\mathbf{p}) = B \log_2(1 + \gamma_k(\mathbf{p})) = q_k^+(\mathbf{p}) - q_k^-(\mathbf{p}) \qquad (1)$$

*with $q_k^+(\mathbf{p}), q_k^-(\mathbf{p}) : \mathbb{R}_+^K \to \mathbb{R}_+$.*

Additional specific assumptions on $q_k^+$ and $q_k^-$ will be introduced in Sections IV and V, when discussing the monotonic fractional programming and the sequential fractional programming frameworks, respectively. For now, it should be stressed that Assumption 1 is very general, and holds at least in the following three notable cases.

(*i*) The typical SINR expression in interference networks takes the general form

$$\gamma_k(\mathbf{p}) = \frac{p_k \alpha_k}{\sigma^2 + \sum\limits_{i=1, i\neq k}^{K} p_i \beta_{i,k}} \qquad (2)$$

where $\sigma^2$ is the power [in W] of the receiver noise (over the bandwidth $B$), $\alpha_k$ is the channel gain over link $k$, whereas $\{\beta_{i,k}\}$ account for the multi-user interference and depend on the other links' channel coefficients as well as on global system parameters. The particular expression of coefficients $\{\alpha_k, \{\beta_{i,k}\}\}$ is determined by the specific system under consideration. From (2), it follows that $q_k^+$ and $q_k^-$ take the form:

$$q_k^+(\mathbf{p}) = \log_2\left(\sigma^2 + p_k\alpha_k + \sum_{i=1, i\neq k}^{K} p_i \beta_{i,k}\right) \qquad (3)$$

$$q_k^-(\mathbf{p}) = \log_2\left(\sigma^2 + \sum_{i=1, i\neq k}^{K} p_i \beta_{i,k}\right). \qquad (4)$$

(*ii*) A considerable extension of (2) is obtained by considering also a self-interference term in the denominator, proportional to the useful power:

$$\gamma_k(\mathbf{p}) = \frac{p_k \alpha_k}{\sigma^2 + p_k\phi_k + \sum\limits_{i=1, i\neq k}^{K} p_i \beta_{i,k}} \qquad (5)$$

with the coefficients $\{\phi_k\}$ also depending on propagation channels and system parameters. A non-zero coefficient $\phi_k$ arises in several relevant instances of communication systems, such as hardware-impaired networks, receivers with imperfect channel state information (CSI), relay-assisted communications, and systems affected by inter-symbol interference. A detailed discussion on the communication scenarios in which the SINR may take the form in (5) can be found in [19]. Given (5), the functions $q_k^+$ and $q_k^-$ are easily found to be:

$$q_k^+(\mathbf{p}) = \log_2\left(\sigma^2 + p_k(\alpha_k + \phi_k) + \sum_{i=1, i\neq k}^{K} p_i \beta_{i,k}\right) \qquad (6)$$

$$q_k^-(\mathbf{p}) = \log_2\left(\sigma^2 + p_k\phi_k + \sum_{i=1, i\neq k}^{K} p_i \beta_{i,k}\right). \qquad (7)$$

(*iii*) A third notable SINR expression is that obtained in vector channels, when linear minimum mean square error (LMMSE) reception is used at the receiver:

$$\gamma_k(\mathbf{p}) = p_k \mathbf{v}_{kk}^H \left(\sigma^2 \mathbf{I}_r + p_k \mathbf{u}_k \mathbf{u}_k^H + \sum_{i=1, i\neq k}^{K} p_i \mathbf{v}_{ki} \mathbf{v}_{ki}^H\right)^{-1} \mathbf{v}_{kk} \qquad (8)$$

where $r$ denotes the dimension of the received signal, $\mathbf{v}_{ki}$ is the $r \times 1$ channel vector between transmitter $i$ and receiver $k$, and the $r \times 1$ vector $\mathbf{u}_k$ accounts for self-interference terms (due to the same reasons as for the SINR expression in (5)). It should be stressed that since LMMSE is the SINR-maximizing linear receive structure [29], and since the receiver choice is clearly decoupled among the different links, LMMSE turns out to be the optimal linear receive structure from an energy-efficient perspective, too. Thus, optimizing the transmit powers assuming the SINR expression (8), means finding the optimal transmit powers after linear receiver design has been already performed and the corresponding optimal SINR has been plugged into the energy efficiency functions. In other words, considering the SINR (8) allows performing joint power and receiver design.

Given (8), the functions $q_k^+$ and $q_k^-$ are expressed as

$$q_k^+(\mathbf{p}) = \log_2\left|\sigma^2\mathbf{I}_r + p_k(\mathbf{v}_{kk}\mathbf{v}_{kk}^H + \mathbf{u}_k\mathbf{u}_k^H) + \sum_{i=1, i\neq k}^{K} p_i\mathbf{v}_{ki}\mathbf{v}_{ki}^H\right| \tag{9}$$

$$q_k^-(\mathbf{p}) = \log_2\left|\sigma^2\mathbf{I}_r + p_k\mathbf{u}_k\mathbf{u}_k^H + \sum_{i=1, i\neq k}^{K} p_i\mathbf{v}_{ki}\mathbf{v}_{ki}^H\right|. \tag{10}$$

In the considered interference network scenario, the EE (measured in bit/Joule) of link $k$ is defined as the ratio of the achievable rate and the total power consumption

$$\mathrm{EE}_k(\mathbf{p}) = \frac{B\log_2(1+\gamma_k(\mathbf{p}))}{\mu_k p_k + \Psi_k} \tag{11}$$

wherein $\mu_k \geq 1$ is the inverse of the power amplifier efficiency of transmitter node $k$ and $\Psi_k$ is the circuit power required to operate link $k$, accounting for the dissipation in analog hardware, digital signal processing, backhaul signaling, and other overhead costs (such as cooling and power supply losses) [5], [30], [31]. Clearly, (11) is a link-centric (or user-centric) performance metric. A network-centric definition of EE must combine the individual energy efficiencies of the different links. Although different approaches have been proposed in the literature, a single definition that unarguably best represents the EE of the whole network is not available, since the different EEs are typically conflicting objectives [5], [32]. Two of the most well-established metrics to measure the network EE are the GEE defined as

$$\mathrm{GEE}(\mathbf{p}) = \frac{\sum_{k=1}^{K} B\log_2(1+\gamma_k(\mathbf{p}))}{\sum_{k=1}^{K} \mu_k p_k + \Psi_k} \tag{12}$$

and the weighted minimum energy efficiency (WMEE) given by

$$\mathrm{WMEE}(\mathbf{p}) = \min_{k=1,\ldots,K} w_k \frac{B\log_2(1+\gamma_k(\mathbf{p}))}{\mu_k p_k + \Psi_k} \tag{13}$$

where the coefficients $w_k \in \mathbb{R}_+$ are used to weigh the EEs of the individual links. Two other possible metrics are the weighted sum energy efficiency (WSEE), defined as

$$\mathrm{WSEE}(\mathbf{p}) = \sum_{k=1}^{K} w_k \frac{B\log_2(1+\gamma_k(\mathbf{p}))}{\mu_k p_k + \Psi_k} \tag{14}$$

and the weighted product energy efficiency (WPEE), defined as

$$\mathrm{WPEE}(\mathbf{p}) = \prod_{k=1}^{K} \left[\frac{B\log_2(1+\gamma_k(\mathbf{p}))}{\mu_k p_k + \Psi_k}\right]^{w_k}. \tag{15}$$

The GEE is the metric with the strongest physical interpretation, as it represents the benefit-cost ratio of the entire network, in terms of global amount of reliably transmitted data and global amount of consumed energy. However, it does not depend on the individual EEs, and therefore does not allow tuning the EE of the individual links according to specific needs. Instead, the WMEE, WSEE, and WPEE are more connected to a multi-objective approach, in which the objectives are the individual EEs [5], [32], which turns out to be quite useful in heterogeneous networks. By suitably choosing the weights, it is possible to prioritize the links that require higher EE, choosing different operating points in the system energy-efficient Pareto region, defined as the region $\mathcal{E}$ containing all feasible $K \times 1$ vectors of the users' EEs:

$$\mathcal{E} = \left\{[\mathrm{EE}_1(\mathbf{p}), \ldots, \mathrm{EE}_K(\mathbf{p})]^T \;:\; \mathbf{p} \in \mathcal{P}\right\}. \tag{16}$$

A known result from multi-objective optimization theory ensures that, by globally maximizing the WMEE for different choices of the weights, it is possible to describe the complete Pareto boundary of the Pareto region (16). In general, this is not possible by maximizing the WSEE or WPEE; for example, the WSEE maximization allows only to describe the convex hull of the Pareto region, and therefore the boundary of (16) can not be completely characterized unless the region is convex. In light of these considerations, the main focus of this article is on the maximization of GEE and WMEE. Nevertheless, it will be shown that most of the techniques developed in the sequel are general enough to apply also to the maximization of WSEE and WPEE.

Given this background, the problem to be solved can be mathematically stated as:

$$\underset{\mathbf{p}}{\text{maximize}} \quad u(\mathbf{p}) \quad \text{s.t.} \quad \mathbf{p} \in \mathcal{P} \tag{17}$$

wherein the objective $u(\mathbf{p})$ is chosen as either the GEE or WMEE (given by (12) or (13)), whereas $\mathcal{P}$ represents the feasible set of the problem given by

$$\mathcal{P} = \left\{\mathbf{p} \in \mathbb{R}_+^K;\; p_k \leq P_{\max,k}, c_k(\mathbf{p}) \geq 0 \;\; \forall k \in \{1, \ldots, K\}\right\} \tag{18}$$

where $c_k(\mathbf{p}) : \mathbb{R}_+^K \to \mathbb{R}_+$ accounts possible additional constraint functions. Similarly to what has been done for SINRs, no particular expression is assumed here for $c_k(\mathbf{p})$, and only the following assumption is made:

**Assumption 2.** *The function $c_k(\mathbf{p})$ can be expressed $\forall k$ as the difference of two non-negative functions, namely:*

$$c_k(\mathbf{p}) = c_k^+(\mathbf{p}) - c_k^-(\mathbf{p}) \tag{19}$$

*with* $c_k^+(\mathbf{p}), c_k^-(\mathbf{p}) : \mathbb{R}_+^K \to \mathbb{R}_+$.

Additional specific assumptions on $c_k^+$ and $c_k^-$ will be introduced in Sections IV and V, when discussing the monotonic fractional programming and the sequential fractional programming frameworks, respectively. Also, observe that the per-user power constraint in (18) well models uplink transmissions. However, the constraint functions $c_k(\mathbf{p})$ can be defined to also enforce a total power constraint $\sum_{k=1}^K p_k \leq P_{\max}$, which is particularly relevant in downlink channels. In addition, the constraint functions $c_k(\mathbf{p})$ naturally model minimum rate constraints and/or interference temperature constraints. These cases are obtained by defining $c_k(\mathbf{p})$ to be the achievable rate of user $k$, or the achievable rate of another network node, whose performance must be protected from the interference from user $k$. It should be observed that while a single constraint function $c_k(\mathbf{p})$ is considered for each network node in the definition of the feasible set of (18), this is due only to notational convenience. The frameworks to be developed immediately apply to scenarios in which multiple constraint functions are enforced for each network node, so that minimum rate constraint and interference temperature constraints can be enforced at the same time on each network node $k$, together with power constraints. Moreover, it is also possible to include additional, network-wide constraint functions, to require a minimum network sum-rate or to limit the maximum interference that the network causes to neighboring communication systems.

A further point to be made is about CSI assumptions. Again, the problem formulation in (17) is very general and does not require any particular CSI assumption. The SINRs $\{\gamma_k\}_{k=1}^K$ can be either based on perfect CSI, or on statistical CSI, or on imperfect CSI. When only statistical CSI is available it is possible to either use the channel second-order statistics in the SINRs expressions, or to replace the rate functions in the numerators of the EE metrics with their ergodic counterparts, which is a necessary approach in fast-fading scenarios. Instead, when imperfect CSI is available due to channel estimation errors or limited feedback, it is possible to consider worst-case rate functions, or to consider the average of the numerators of the EE with respect to the channel estimation error distribution. The approaches to be developed in the sequel can work in these scenarios, too, and only require very mild assumptions which will be detailed later.

Regardless of the choice of the EE metric $u(\mathbf{p})$ and of the particular constraint functions to be enforced, the optimization problem in (17) belongs to the class of fractional programming problems [5]. Such problems can be solved with polynomial complexity only if the numerator and denominator of the fraction to maximize are respectively concave and convex, and if the feasible set is also convex [5]. Unfortunately, this requirement is not fulfilled in interference-limited networks as it follows from the general SINR expressions given above. In fact, the functions $q_k^-$ are always non-zero whenever multi-user interference is present, and this causes the links' achievable rates (i.e., the numerators of the individual EEs) to be non-concave functions of $\mathbf{p}$. As a result, fractional programs are in general NP-hard in interference-limited scenarios [4], [5], which calls for new tools to complement and extend the potentialities of fractional programming theory. In Section IV, we aim at solving (17) at the expense of computational complexity by combining fractional programming with monotonic optimization. Then, we look for local solutions of (17) that can be obtained with affordable complexity. This is accomplished in Section V by merging fractional programming with sequential optimization. Before turning to the development of these new optimization frameworks, the next section will provide the necessary mathematical preliminaries and definitions.

## III. Mathematical Preliminaries

This section provides a background on the optimization theories to be used in the remainder. Section III-A gives a short review of fractional programming theory [33] and also provides some background (see [34] for more details) to understand the complexity arguments mentioned at the end of Section II. Then, Section III-B gives an overview of monotonic optimization [21], [22]. Finally, Section III-C briefly discusses the framework of sequential optimization.

### A. Fractional programming

For a more comprehensive overview of fractional programming for EE maximation, the reader is referred to [5].

**Definition 1** (Generalized fractional program)**.** *Let $\mathcal{D} \subseteq \mathbb{R}^N$ and consider the functions $f_k : \mathcal{D} \to \mathbb{R}$ and $g_k : \mathcal{D} \to \mathbb{R}_{++}$, with $k = 1, \ldots, K$. A generalized fractional program is the optimization problem defined as*

$$\underset{\mathbf{x}}{\text{maximize}} \min_{k=1,\ldots,K} \frac{f_k(\mathbf{x})}{g_k(\mathbf{x})} \quad \text{s.t.} \quad \mathbf{x} \in \mathcal{D}. \tag{20}$$

*If $K = 1$, then the above problem reduces to the so-called single-ratio fractional program:*

$$\underset{\mathbf{x}}{\text{maximize}} \frac{f_1(\mathbf{x})}{g_1(\mathbf{x})} \quad \text{s.t.} \quad \mathbf{x} \in \mathcal{D}. \tag{21}$$

Since the objective function in (20) is in general not concave, standard convex optimization algorithms are not guaranteed to solve (20) and specific algorithms are required. Towards this end, we have the following main result.

**Proposition 1.** *[35], [36]. A vector $\mathbf{x}^\star \in \mathcal{D}$ solves* (20) *if and only if*

$$\mathbf{x}^\star = \arg\max_{\mathbf{x}\in\mathcal{D}} \left\{ \min_{k=1,\ldots,K} \left[ f_k(\mathbf{x}) - \lambda^\star g_k(\mathbf{x}) \right] \right\} \tag{22}$$

*with $\lambda^\star$ being the unique zero of the auxiliary function $F(\lambda)$:*

$$F(\lambda) = \max_{\mathbf{x}\in\mathcal{D}} \min_{k=1,\ldots,K} \{ f_k(\mathbf{x}) - \lambda g_k(\mathbf{x}) \}. \tag{23}$$

This result allows one to solve (20) by finding the unique zero of $F(\lambda)$. To this end, the most widely used algorithm is the (Generalized, if $K > 1$) Dinkelbach's algorithm [34], [36], reported in Algorithm 1.

A critical point about Algorithm 1 is that it converges to the global optimum of the corresponding instance of the fractional problem only provided that (22) can be *globally* solved at each iteration. If $f$ is concave, $g$ is convex, and all

**Algorithm 1** Generalized Dinkelbach's algorithm

Initialize $\lambda_0$ with $F(\lambda_0) \geq 0$, $j = 0$;
**while** $F(\lambda_j) > \varepsilon$ **do**
  Solve the problem:
  $$\mathbf{x}_j^\star = \arg\max_{\mathbf{x}\in\mathcal{D}} \left\{ \min_{k=1,\ldots,K} \left[ f_k(\mathbf{x}) - \lambda_j g_k(\mathbf{x}) \right] \right\};$$
  $\lambda_{j+1} = \min_{k=1,\ldots,K} \frac{f_k(\mathbf{x}_j^\star)}{g_k(\mathbf{x}_j^\star)}$;
  $j = j + 1$;
**end while**

constraints are also convex, then this can be accomplished with polynomial-time complexity, since each subproblem is a concave maximization subject to convex constraints.[1] Moreover, Algorithm 1 exhibits a super-linear convergence rate, since the update rule for $\lambda$ follows Newton's method applied to the function $F(\lambda)$ [34].

Instead, if (22) is not a non-convex problem, the complexity of Algorithm 1 becomes significantly higher, because globally solving (22) can not be handled by the well-developed theory of convex optimization, but the use of global optimization algorithms is required. However, standard global optimization methods operate by exploring the whole feasible set [39], with a prohibitive computational complexity, even for small problem instances, and with a convergence that is only guaranteed if the functions have a limited variability (e.g., Lipschitz continuity [27]).

Finally, it should be observed that other approaches exist to tackle a fractional problem, but, similarly to Dinkelbach's method, all of them require the concavity of $f$, the convexity of $g$, plus the convexity of the feasible set, in order to require only the solution of convex problems. In particular, we mention the Charnes-Cooper transform method, which is able to convert a fractional problem into an equivalent convex problem, provided $f$ is concave, $g$ is convex, and the feasible set is also convex [40]. Unlike Dinkelbach's algorithm, this approach is not iterative, since a single convex problem must be solved, but this equivalent problem has one additional variable and constraint compared to Problem (22). Moreover, the Charnes-Cooper transform involves the use of the perspective function of the numerator $f$, which is not suited to the global optimization algorithm to be developed in Section IV.

### B. Monotonic optimization

Monotonic optimization is a relatively recent global optimization framework, which exploits monotonicity or hidden monotonicity structures in the objective and constraints to reduce computational complexity and provide a guaranteed convergence [21], [22]. The basic idea is that if the objective to be maximized is increasing in all optimization variables, then it is not necessary to explore the complete feasible set of the problem, but only its outer boundary. This concept is made formal in the rest of this section.

[1]Such a class of problems can be solved with polynomial-time complexity [37], [38]

**Definition 2** (Monotonicity in $\mathbb{R}^N$)**.** *A function $f : \mathbb{R}^N \to \mathbb{R}$ is monotonically increasing if $f(\mathbf{y}) \geq f(\mathbf{x})$ when $\mathbf{y} \succeq \mathbf{x}$.*

**Definition 3** (Hyper-rectangle in $\mathbb{R}^N$)**.** *Let $\mathbf{a},\ \mathbf{b} \in \mathbb{R}^N$ with $\mathbf{a} \preceq \mathbf{b}$. Then, the set of all $\mathbf{x} \in \mathbb{R}^N$ such that $\mathbf{a} \preceq \mathbf{x} \preceq \mathbf{b}$ is a hyper-rectangle in $\mathbb{R}^N$ and is denoted by $[\mathbf{a}, \mathbf{b}]$.*

**Definition 4** (Normal and Co-normal sets)**.** *A set $\mathcal{S} \subset \mathbb{R}^N$ is normal if $\forall \mathbf{x} \in \mathcal{S}$, the hyper-rectangle $[\mathbf{0}, \mathbf{x}]$ belongs to $\mathcal{S}$. A set $\mathcal{S}_c \subset \mathbb{R}^N$ is co-normal in $[\mathbf{0}, \mathbf{b}]$ if $\forall \mathbf{x} \in \mathcal{S}_c$, then $[\mathbf{x}, \mathbf{b}] \subset \mathcal{S}_c$.*

A given function $h : \mathbb{R}^N \to \mathbb{R}$ defines a normal or a co-normal set if the following results hold true:

**Proposition 2.** *[21] The set $\mathcal{S} = \{\mathbf{x} \in \mathbb{R}^N : h(\mathbf{x}) \leq 0\}$ is normal and closed if $h$ is lower semi-continuous and increasing. The set $\mathcal{S}_c = \{\mathbf{x} \in \mathbb{R}^N : h(\mathbf{x}) \geq 0\}$ is co-normal and closed if $h$ is upper semi-continuous and increasing.*

**Definition 5** (Monotonic optimization)**.** *A monotonic optimization problem in canonical form is defined as*

$$\underset{\mathbf{x}}{\text{maximize}}\ f(\mathbf{x}) \quad \text{s.t.} \quad \mathbf{x} \in \mathcal{S} \cap \mathcal{S}_c \tag{24}$$

*where $f : \mathbb{R}^N \to \mathbb{R}$ is an increasing function, $\mathcal{S} \subset [\mathbf{0}, \mathbf{b}]$ is a compact, normal set with nonempty interior, and $\mathcal{S}_c$ is a closed co-normal set in $[\mathbf{0}, \mathbf{b}]$.*

The main result of monotonic optimization theory states that the solution to (24) lies on the upper boundary of $\mathcal{S} \cap \mathcal{S}_c$ [21, Proposition 7]. Therefore, methods like the polyblock algorithm [21] and the BRB algorithm [22] can be used to globally solve (24) by searching only on the upper boundary of the feasible set, thus drastically simplifying the problem. Nevertheless, we remark that the complexity of monotonic optimization methods is still exponential in the number of variables. However, as already observed, it is much lower than general global optimization methods, which do not exploit any monotonicity structure [21]. This makes monotonic optimization attractive for the development of an off-line framework to benchmark any suboptimal method for solving (24).

More in detail, the complexity of monotonic optimization methods depends on the number of iterations required to search the boundary of the feasible set, and on the complexity to evaluate the objective function on each given point on the boundary of the feasible set. The former in general grows exponentially with the number of variables, whereas the latter is difficult to evaluate in general, since it depends on the particular functional form of the objective function. Nevertheless, it must be emphasized again that using standard global optimization algorithms would be significantly more complex, since not just the frontier, but the whole feasible set should be searched.

### C. Sequential optimization

Sequential optimization is a powerful tool that provides the means to generate candidate solutions of non-convex optimization problems with affordable complexity [41], while at the same time satisfying theoretical optimality claims. This statement is made precise in the following result, which readily follows from [41]:

**Proposition 3.** *Let $\mathcal{F}$ be a maximization problem with differentiable objective $f_0(\mathbf{x})$ and constraints $f_i(\mathbf{x}) \geq 0$ $\forall i \in \{1,\ldots,I\}$, and with a compact feasible set. Let $\mathcal{G}_j$ be a maximization problem with differentiable objective $g_{0,j}(\mathbf{x})$ and constraints $g_{i,j}(\mathbf{x}) \geq 0$, $\forall i \in \{1,\ldots,I\}$, with compact feasible set, and optimal solution $\mathbf{x}_j^\star$. Assume that $\forall j$ and $\forall i \in \{1,\ldots,I\}$ $g_{i,j}(\cdot)$ satisfies the following two properties:*

1) $g_{i,j}(\mathbf{x}) \leq f_i(\mathbf{x})$ $\forall \mathbf{x}$;
2) $g_{i,j}(\mathbf{x}_{j-1}^\star) = f_i(\mathbf{x}_{j-1}^\star)$.

*Then, the sequence $\{f_0(\mathbf{x}_j^\star)\}$ is monotonically increasing and converges to a finite limit $g$. Next, assume the following third property is also satisfied $\forall j$ and $\forall i \in \{1,\ldots,I\}$:*

3) $\nabla g_{i,j}(\mathbf{x}_{j-1}^\star) = \nabla f_i(\mathbf{x}_{j-1}^\star)$ .

*Then, under suitable constraint qualifications, every limit point of $\{\mathbf{x}\}_\mathbf{j}$ that achieves the objective value $g$ fulfills the KKT conditions of the original problem $\mathcal{F}$.*

Proposition 3 shows that by solving the sequence of approximate problems $\{\mathcal{G}_j\}$, one can generate a sequence of feasible points $\mathbf{x}_j^\star$ that monotonically increases the value of the original objective $f_0$. Moreover, the limit value $g$ is the value that the original objective attains at a KKT point of $\mathcal{F}$. The critical issue for this tool to be of practical use, is to find suitable approximate problems $\{\mathcal{G}_j\}$ fulfilling the assumptions of Proposition 3, while at the same being easier to solve than the original problem. Notable cases in which this practical requirement holds are those in which the approximate problem is a concave or pseudo-concave maximization subject to convex constraints. Instead, as for the number of iterations required for the sequential method to converge, no general formulas are available, since this depends on the particular structure of the problem.

It should also be mentioned that recent works provided notable extensions of the seminal result from [41]. In [42] it is shown that not only the sequence $\{f_0(\mathbf{x}_j^\star)\}$ converges, but also that every limit point of the sequence $\{\mathbf{x}\}_\mathbf{j}^\star$ fulfills the KKT conditions of the original problem, provided that the original objective function $f_0$ is strictly concave (assuming the case of maximization problems) and under regularity conditions. Instead, [43] shows that every limit point of the sequence $\{\mathbf{x_j}\}_\mathbf{j}^\star$ is a stationary point of the original problem also when the original objective is not strictly concave, provided the original problem has a convex feasible set, and only the objective function is approximated.

## IV. Global optimality: Monotonic Fractional Programming

As mentioned above, among global optimization algorithms, monotonic optimization provides attractive complexity and convergence properties. However, it can not be directly employed to solve (17), because the EEs are not monotone functions of $\mathbf{p}$ in the sense of Definition 2. This section shows how this difficulty can be overcome by an interplay of fractional programming and monotonic optimization, provided that the following assumption holds:

**Assumption 3.** *The functions $q_k^+(\mathbf{p})$ and $q_k^-(\mathbf{p})$ in* (1)*, and the functions $c_k^+(\mathbf{p})$ and $c_k^-(\mathbf{p})$ in* (19) *are monotonic functions $\forall k \in \{1,\ldots,K\}$ as stated in Definition 2.*

In other words, it is assumed that all achievable rates and constraint functions can be written as the difference of monotonic functions. Observe that no assumption on the concavity or convexity of $q_k^+(\mathbf{p})$, $q_k^-(\mathbf{p})$, $c_k^+(\mathbf{p})$, and $c_k^-(\mathbf{p})$ is made.

### A. GEE maximization

GEE maximization belongs to the class of single-ratio fractional problems. Thus, finding its solution by Dinkelbach's algorithm requires to solve the following auxiliary problem at iteration $j$:

$$\underset{\mathbf{p}}{\text{maximize}} \sum_{k=1}^{K} B\log_2(1+\gamma_k) - \lambda_j(\mu_k p_k + \Psi_k) \quad \text{s.t. } \mathbf{p} \in \mathcal{P} \tag{25}$$

for a given positive $\lambda_j$. Note that, at a first sight, the above problem is not a monotonic optimization problem in canonical form because:

- The objective function is not monotonic, since the achievable rates $\log_2(1+\gamma_k)$ are not increasing functions of $\mathbf{p}$, and since the negative term is in fact decreasing in the interfering powers.
- The constraint set is not guaranteed to be the intersection of a normal and a co-normal set, since the difference of two increasing functions is in general not increasing.

However, (25) exhibits a hidden monotonicity structure as shown in the following proposition:

**Proposition 4.** *If Assumption 3 holds true, then* (25) *can be expressed as a monotonic optimization problem in canonical form.*

*Proof:* Observe that (25) can be equivalently written as

$$\underset{\mathbf{p}}{\text{maximize}} \quad q^+(\mathbf{p}) - q^-(\mathbf{p},\lambda_j) \quad \text{s.t.} \quad \mathbf{p} \in \mathcal{P} \tag{26}$$

wherein $q^+(\mathbf{p})$ and $q^-(\mathbf{p},\lambda_j)$ are increasing in $\mathbf{p}$ and given by

$$q^+(\mathbf{p}) = \sum_{k=1}^{K} q_k^+(\mathbf{p}) \tag{27}$$

$$q^-(\mathbf{p},\lambda_j) = \sum_{k=1}^{K} q_k^-(\mathbf{p}) + \lambda_j\big(\mu_k p_k + \Psi_k\big). \tag{28}$$

Next, define $\mathbf{p}_{\max} = [P_{\max,1},\ldots,P_{\max,K}]$ and introduce the auxiliary variable $t = q^-(\mathbf{p}_{\max},\lambda_j) - q^-(\mathbf{p},\lambda_j)$. Then, for any given $\lambda_j$, (26) can be rewritten as

$$\underset{(t,\,\mathbf{p})}{\text{maximize}} \quad q^+(\mathbf{p}) + t \tag{29}$$

$$\text{s.t.} \quad (t,\,\mathbf{p}) \in \mathcal{P} \cap \mathcal{Q} \tag{30}$$

with

$$\mathcal{Q} = \left\{ (t,\mathbf{p}) : \begin{array}{l} 0 \leq t + q^-(\mathbf{p},\lambda_j) \leq q^-(\mathbf{p}_{\max},\lambda_j) \\ 0 \leq t \leq q^-(\mathbf{p}_{\max},\lambda_j) - q^-(\mathbf{0}_K,\lambda_j) \end{array} \right\}.$$

Problem (29) is not a monotonic problem yet, because the constraint functions $c_k(\mathbf{p})$ are expressed as the difference of increasing functions. To overcome this problem, observe that the set of constraints $c_k(\mathbf{p}) \geq 0$ with $k = 1, \ldots, K$, can be equivalently rewritten as the following single constraint:

$$\min_{k=1,\ldots,K} \left[c_k^+(\mathbf{p}) - c_k^-(\mathbf{p})\right] \geq 0 \iff \tag{31}$$

$$\min_{k=1,\ldots,K} \left[c_k^+(\mathbf{p}) - \left(\sum_{i=1}^{K} c_i^-(\mathbf{p}) - \sum_{i=1, i\neq k}^{K} c_i^-(\mathbf{p})\right)\right] = \tag{32}$$

$$\underbrace{\min_{k=1,\ldots,K} \left[c_k^+(\mathbf{p}) + \sum_{i=1, i\neq k}^{K} c_i^-(\mathbf{p})\right]}_{c^+(\mathbf{p})} - \underbrace{\sum_{i=1}^{K} c_i^-(\mathbf{p})}_{c^-(\mathbf{p})} \geq 0 \tag{33}$$

which is the difference of the two increasing functions $c^+(\mathbf{p})$ and $c^-(\mathbf{p})$. Similarly as above, we can thus introduce the auxiliary variable $s$ and reformulate (29) as

$$\begin{aligned} \underset{(t,s,\mathbf{p})}{\text{maximize}} \quad & q^+(\mathbf{p}) + t \\ \text{s.t.} \ & (t,\, \mathbf{p}) \in \mathcal{Q}\,,\ 0 \leq s \leq c^-(\mathbf{p}_{\max}) - c^-(\mathbf{0}_K) \\ & c^-(\mathbf{p}) + s \leq c^-(\mathbf{p}_{\max})\,, c^+(\mathbf{p}) + s \geq c^-(\mathbf{p}_{\max}). \end{aligned} \tag{34}$$

In order to complete the proof, it remains to verify that (34) fulfills Definition 5, thus being a monotonic problem in canonical form. To this end, let us first observe that the objective of (34) is monotonic in $(t, s, \mathbf{p})$.

Next, to show that the feasible set of (34) is the intersection of a normal and a co-normal set, let us observe that, for any $\mathbf{p}$ in the feasible set, we have that

$$q^-(\mathbf{0}_K, \lambda_j) \leq q^-(\mathbf{p}, \lambda_j) \tag{35}$$

$$c^-(\mathbf{0}_K) \leq c^-(\mathbf{p}). \tag{36}$$

As a consequence, the feasible set of (34) can be written as the intersection of the following two sets:

$$\mathcal{S} = \Big\{(t,s,\mathbf{p}) : \mathbf{p} \preceq \mathbf{p}_{\max}, t + q^-(\mathbf{p},\lambda_j) \leq q^-(\mathbf{p}_{\max},\lambda_j), \\ s + c^-(\mathbf{p}) \leq c^-(\mathbf{p}_{\max})\Big\} \tag{37}$$

$$\mathcal{S}_c = \left\{(t,s,\mathbf{p}) : \mathbf{p} \succeq \mathbf{0}_K, t \geq 0, s + c^+(\mathbf{p}) \geq c^-(\mathbf{p}_{\max})\right\}\ . \tag{38}$$

Then, since all the constraint functions in (37) and (38) are monotonic and continuous, by virtue of Proposition 2 and employing again (35), it follows that $\mathcal{S}$ and $\mathcal{S}_c$ are normal and co-normal sets in the hyper-rectangle given by:

$$[0, q^-(\mathbf{p}_{\max},\lambda_j) - q^-(\mathbf{0}_K,\lambda_j)] \times [c^-(\mathbf{p}_{\max}) - c^-(\mathbf{0}_K)] \\ \times [\mathbf{0}_K, \mathbf{p}_{\max}]. \tag{39}$$

This completes the proof. ■

### *B. WMEE maximization*

WMEE maximization belongs to the class of generalized fractional programs and requires to solve the following auxiliary problem at iteration $j$:

$$\begin{aligned} \underset{\mathbf{p}}{\text{maximize}} \min_{k=1,\ldots,K} \ & q_k^+(\mathbf{p}) - q_k^-(\mathbf{p}) - \lambda_j(\mu_k p_k + \Psi_k) \\ \text{s.t.} \quad & \mathbf{p} \in \mathcal{P}. \end{aligned} \tag{40}$$

As in the case of GEE maximization, the objective function is not monotonic. However, the following result can be proved:

**Proposition 5.** *If Assumption 3 holds true, then* (40) *can be expressed as a monotonic problem in canonical form.*

*Proof:* Let $\nu_k(\mathbf{p}, \lambda_j) = q_k^-(\mathbf{p}) + \lambda_j(\mu_k p_k + \Psi_k)$ such that we may rewrite the objective function as

$$\begin{aligned} & q_k^+(\mathbf{p}) - q_k^-(\mathbf{p}) - \lambda_j(\mu_k p_k + \Psi_k) \\ & = q_k^+(\mathbf{p}) - \left(\sum_{i=1}^{K} \nu_i(\mathbf{p},\lambda_j) - \sum_{i=1, i\neq k}^{K} \nu_i(\mathbf{p},\lambda_j)\right) \\ & = \left(q_k^+(\mathbf{p}) + \sum_{i=1, i\neq k}^{K} \nu_i(\mathbf{p},\lambda_j)\right) - \sum_{i=1}^{K} \nu_i(\mathbf{p},\lambda_j). \end{aligned} \tag{41}$$

Then, introduce $t = \sum_{i=1}^{K} \nu_i(\mathbf{p}_{\max}, \lambda_j) - \sum_{i=1}^{K} \nu_i(\mathbf{p}, \lambda_j)$ and reformulate (40) as

$$\begin{aligned} \underset{(t,\mathbf{p})}{\text{maximize}} \min_{k=1,\ldots,K} \ & q_k^+(\mathbf{p}) + \sum_{i=1, i\neq k}^{K} \nu_i(\mathbf{p},\lambda_j) + t \\ \text{s.t.} \quad & (t,\, \mathbf{p}) \in \mathcal{P} \cap \mathcal{Q}' \end{aligned} \tag{42}$$

with

$$\mathcal{Q}' = \left\{(t,\mathbf{p}) : \begin{array}{l} 0 \leq t \leq \sum_{i=1}^{K} \nu_i(\mathbf{p}_{\max},\lambda_j) - \nu_i(\mathbf{p},\lambda_j) \\ 0 \leq t \leq \sum_{i=1}^{K} \nu_i(\mathbf{p}_{\max},\lambda_j) - \nu_i(\mathbf{0}_K,\lambda_j) \end{array}\right\}.$$

Reformulating the constraints $c_k(\mathbf{p}) \geq 0\ \forall k$ as the single constraint (33), (42) is shown to be a monotonic problem in canonical form by using the same arguments adopted in the proof of Proposition 4. ■

Observe that the results of Propositions 4 and 5 do not specifically require an affine power consumption model, as the one adopted in this article. Indeed, the above results can be extended to any power consumption model such that $q^-(\mathbf{p}, \lambda_j)$ (for GEE maximization) and $\nu_k(\mathbf{p}, \lambda_j)$ (for WMEE maximization) are monotonic in $\mathbf{p}$. Moreover, the above result applies also to scenarios in which each user has $N_c > 1$ constraint functions $\{c_{k,i}(\mathbf{p})\}_{i=1}^{N_c}$, provided $c_{k,i}$ can be still expressed as the difference of two monotonic functions.

### *C. WSEE and WPEE maximization*

We now look at the maximization of WSEE and WPEE, which belong to the classes of sum of ratios and product of ratios problems, respectively. Both are hard to solve even if all numerators are concave, all denominators are convex, and the feasible set is convex [5], [44]. Nevertheless, the proposed monotonic fractional programming framework can also be used to solve WSEE and WPEE maximization problems[2]. To

[2]Alternative approaches for the maximization of sum-of-ratios problems have recently appeared in [45], [46]

begin with, observe that the WSEE can be expressed as a single ratio:

$$\text{WSEE}(\mathbf{p}) = \sum_{k=1}^{K} \frac{q_k^+(\mathbf{p}) - q_k^-(\mathbf{p})}{\mu_k p_k + \Psi_k} \tag{43}$$

$$= \frac{\sum_{k=1}^{K} \left(q_k^+(\mathbf{p}) - q_k(\mathbf{p})^-\right) \prod_{i \neq k}(\mu_i p_i + \Psi_i)}{\prod_{k=1}^{K} (\mu_k p_k + \Psi_k)}. \tag{44}$$

Since the product of increasing functions is still an increasing function, (44) turns out to be a single ratio whose numerator is the difference of increasing functions, and the denominator is an increasing function. Hence, the method adopted in Section IV-A can also be used to globally maximize the WSEE.

The same property holds for WPEE maximization since (15) can be reformulated as

$$\text{WPEE}(\mathbf{p}) = \frac{\prod_{k=1}^{K} \left(q_k^+(\mathbf{p}) - q_k^-(\mathbf{p})\right)}{\prod_{k=1}^{K} (\mu_k p_k + \Psi_k)}. \tag{45}$$

Expanding the products in the numerator of the above function, we obtain again the difference of two increasing functions, while the denominator is clearly increasing with respect to $\mathbf{p}$.

**Remark 1.** While this section has focused on power control or on power control and linear receiver design (when the SINR (8) is considered), the approaches described here can be formally extended to transmit beamforming problems, too, upon applying domain changes techniques [25], [26], which enable to reformulate beamforming matrix-variate problems into vector-variate problems formally similar to the power control problems studied in this section, and with a number of variables again equal to the number of links. The resulting problems can be globally solved again by means of the polyblock algorithm. However, the problem feasible set after the domain change technique in general is not available in closed-form, which increases the complexity of the polyblock method.

## V. First-order optimality: Sequential Fractional Programming

As observed in Section IV, despite enjoying a lower complexity than standard global optimization algorithms, the complexity of the proposed monotonic fractional programming framework is still exponential. Motivated by the need of providing also a practical optimization framework for large networks, in this section fractional programming is combined with the sequential optimization results from Proposition 3. This results in a novel sequential fractional programming framework able to compute candidate solutions of EE problems with affordable complexity, while at the same time fulfilling theoretical optimality claims. In particular, the proposed algorithm yields points fulfilling the KKT first-order optimality conditions of the GEE and WMEE maximization problems[3]. Section VI will provide numerical evidence that the sequential fractional programming approach actually attains global optimality, as it finds the same solution as the globally optimal solution computed by means of the monotonic fractional programming framework in Section IV.

This is not the first time that fractional programming and sequential optimization are used together to solve EE optimization problems [17], [19]. In [17], the two theories are employed to maximize the GEE as well as the WPEE under the assumption that the SINR takes the form in (2) whereas (5) is adopted in [19], also including rate constraints in the analysis and addressing WMEE maximization. In these works, the approximate problems $\mathcal{G}_j$ required by Proposition 3 are obtained by means of a logarithmic approximation of the achievable rate plus a change of variable. Compared to [17], [19], a different approach is pursued here that has the following main advantages:

*i*) It is more general since it can handle SINRs in the form of (8), besides those given as (2) and (5), whereas the approaches from [17], [19] do not apply to the SINR expression in (8). It should be stressed that this is not a minor accomplishment, since the SINR in (8) arises when joint power control and receiver design is carried out. Thus, being able to handle the SINR in (8) enables to extend the framework from power control problems, to joint power control and receiver design problems.

*ii*) Another extension compared to available literature is that the proposed framework can include multiple types of constraint functions. Previous works which employ sequential optimization for energy-efficient resource allocation do not enforce interference temperature constraints. This appears a relevant scenario, since it arises for example in device-to-device communications, one strong candidate technology for 5G networks. Instead, as already mentioned, our framework is broad enough to include interference temperature constraints by specializing accordingly the constraint functions $c_k$.

*iii*) As demonstrated in the two points described above, the optimization framework developed here encompasses many relevant scenarios, in terms of both objective functions, and constraint functions, whereas previous works have considered only specific instances of communication systems and related problem formulations. This makes our framework suitable for many important instances of communication systems, including the leading 4G technologies, and the strongest 5G candidate techniques. In fact, in contrast to previous works, we take a higher-level approach, working with generic objectives $q_k$ and constraints $c_k$, and identifying general assumptions that these functions must fulfill in order for the sequential and fractional tools to be used.

In this section, besides the Assumptions in Section II, we only require the following assumption:

**Assumption 4.** *The functions $q_k^+(\mathbf{p})$ and $q_k^-(\mathbf{p})$ in* (1), *and the functions $c_k^+(\mathbf{p})$ and $c_k^-(\mathbf{p})$ in* (19) *are concave functions of $\mathbf{p}$ $\forall k$.*

In other words, it is only required that all the achievable rates and constraint functions can be expressed as the difference of concave functions. It should be observed that

[3]As for the WMEE problem, the KKT conditions refer to the equivalent epigraph-form representation, which removes the non-differentiability of the WMEE function.

Assumption 3 is not required in the sequel. Indeed, in general no relation exists between concavity and monotonicity, since a function can be concave even if it is not monotonic, or vice versa. In this section we only assume the concavity of the functions $q_k^+(\mathbf{p})$ and $q_k^-(\mathbf{p})$, without any claim about their monotonicity. As in Section IV, we consider the GEE and WMEE maximizations separately.

### A. GEE maximization

By virtue of Assumptions 1 and 2, the GEE maximization problem can be cast as:

$$\begin{aligned} \underset{\mathbf{p}}{\text{maximize}} \quad & \frac{\sum_{k=1}^K q_k^+(\mathbf{p}) - q_k^-(\mathbf{p})}{\sum_{k=1}^K \mu_k p_k + \Psi_k} && \text{(46a)} \\ \text{s.t.} \quad & 0 \leq p_k \leq P_{\max,k} \ \forall k && \text{(46b)} \\ & c_k^+(\mathbf{p}) - c_k^-(\mathbf{p}) \geq 0 \ \forall k. && \text{(46c)} \end{aligned}$$

If Assumption 4 holds true, then the numerator of (46a) and the constraint functions in (46c) are expressed as the difference of concave functions, and therefore are not concave in general. As pointed out in Section III-A, this prevents from directly using fractional programming to solve (46). To circumvent this issue, we exploit its hidden structure and obtain the following main result:

**Proposition 6.** *For any given $\mathbf{p}_j$, denote by $\mathcal{G}_j$ the optimization problem*

$$\underset{\mathbf{p}}{\text{maximize}} \ \frac{\sum_{k=1}^K q_k^+(\mathbf{p}) - \left[q_k^-(\mathbf{p}_j) + \left(\nabla_{\mathbf{p}} q_k^-|_{\mathbf{p}=\mathbf{p}_j}\right)^T (\mathbf{p}-\mathbf{p}_j)\right]}{\sum_{k=1}^K \mu_k p_k + \Psi_k} \quad \text{(47a)}$$

$$\text{s.t.} \ 0 \leq p_k \leq P_{\max,k} \ \forall k \quad \text{(47b)}$$

$$c_k^+(\mathbf{p}) - \left[c_k^-(\mathbf{p}_j) + \left(\nabla_{\mathbf{p}} c_k^-|_{\mathbf{p}=\mathbf{p}_j}\right)^T (\mathbf{p}-\mathbf{p}_j)\right] \geq 0 \ \forall k \quad \text{(47c)}$$

*and call $\mathbf{p}_j^\star$ its optimal solution. If $\mathbf{p}_j = \mathbf{p}_{j-1}^\star \ \forall j \geq 1$, and $\mathbf{p}_0$ is any feasible power vector, then $\{\text{GEE}(\mathbf{p}_j^\star)\}_j$ is monotonically increasing and converges to a value $\overline{GEE}$. Moreover, any limit point of the sequence $\{\text{GEE}(\mathbf{p}_j^\star)\}_j$ that achieves $\overline{GEE}$ fulfills the KKT optimality conditions of* (46) *under suitable constraint qualifications*[4].

*Proof:* Recall that any concave function is upper-bounded by its first-order Taylor expansion at any point. Since $q_k^-(\mathbf{p})$ and $c_k^-(\mathbf{p})$ are concave functions, for any power vector $\mathbf{p}_j$ we thus have that

$$\begin{aligned} q_k^+(\mathbf{p}) - q_k^-(\mathbf{p}) &\geq q_k^+(\mathbf{p}) - \left[q_k^-(\mathbf{p}_j) + \left(\nabla_{\mathbf{p}} q_k^-|_{\mathbf{p}=\mathbf{p}_j}\right)^T (\mathbf{p}-\mathbf{p}_j)\right] \\ c_k^+(\mathbf{p}) - c_k^-(\mathbf{p}) &\geq c_k^+(\mathbf{p}) - \left[c_k^-(\mathbf{p}_j) + \left(\nabla_{\mathbf{p}} c_k^-|_{\mathbf{p}=\mathbf{p}_j}\right)^T (\mathbf{p}-\mathbf{p}_j)\right] \end{aligned} \quad \text{(48)}$$

Hence, (47a) and (47c) are lower bounds of (46a) and (46c), respectively. Next, since the lower bounds in (48) are tight when evaluated in $\mathbf{p}_j$, it immediately follows that (47a) and (47c) are equal to (46a) and (46c), respectively, for $\mathbf{p} = \mathbf{p}_j$. Similarly, it can be shown that the gradients of (47a) and (47c) are equal to those of (46a) and (46c), for $\mathbf{p} = \mathbf{p}_j$. Thus, (47) fulfills all the assumptions of Proposition 3 with respect to (46), which completes the proof of this proposition. ■

For any $\mathbf{p}_j$, (47a) has a concave numerator, and affine denominator, while the constraint functions in (47b) and (47c) are all affine or concave. As a result, (47) is a single-ratio problem, which can be globally solved by means of fractional programming theory. In particular, it is possible to either use Algorithm 1 directly solving (47), or Charnes-Cooper transform can be used to further reformulate (47) into a an equivalent convex program [40].

**Remark 2.** *Observe that the gradients in* (48) *can be computed in closed form, once the SINR expression is specified. If the SINRs are expressed as in* (5)*, then we have*

$$\frac{\partial q_k^-}{\partial p_\ell} = \begin{cases} \dfrac{\phi_k}{\ln(2)\left(\sigma^2 + \phi_k p_k + \sum_{i=1, i\neq k}^K \beta_{i,k} p_i\right)}, & \ell = k, \\ \dfrac{\beta_{\ell,k}}{\ln(2)\left(\sigma^2 + \phi_k p_k + \sum_{i=1, i\neq k}^K \beta_{i,k} p_i\right)} & \ell \neq k. \end{cases} \quad \text{(49)}$$

*The formulas for the SINR* (2) *can be obtained by setting $\phi_k = 0$ in* (49)*. On the other hand, if the SINR takes the expression in* (8)*, then we have*[5]

$$\frac{\partial q_k^-}{\partial p_\ell} = \begin{cases} \dfrac{1}{\ln(2)} \mathbf{u}_k^H \left(\sigma^2 \mathbf{I}_r + p_k \mathbf{u}_k \mathbf{u}_k^H + \sum_{i=1, i\neq k}^K p_i \mathbf{v}_{ki} \mathbf{v}_{ki}^H\right)^{-1} \mathbf{u}_k, & \ell = k, \\ \dfrac{1}{\ln(2)} \mathbf{v}_{k\ell}^H \left(\sigma^2 \mathbf{I}_r + p_k \mathbf{u}_k \mathbf{u}_k^H + \sum_{i=1, i\neq k}^K p_i \mathbf{v}_{ki} \mathbf{v}_{ki}^H\right)^{-1} \mathbf{v}_{k\ell}, & \ell \neq k. \end{cases} \quad \text{(50)}$$

### B. WMEE maximization

If Assumptions 1 and 2 are satisfied, the WMEE maximization problem can be cast as

$$\begin{aligned} \underset{\mathbf{p}}{\text{maximize}} \quad & \min_{k=1,\ldots,K} \frac{q_k^+(\mathbf{p}) - q_k^-(\mathbf{p})}{\mu_k p_k + \Psi_k} && \text{(51a)} \\ \text{s.t.} \quad & 0 \leq p_k \leq P_{\max,k} \ \forall k && \text{(51b)} \\ & c_k^+(\mathbf{p}) - c_k^-(\mathbf{p}) \geq 0 \ \forall k. && \text{(51c)} \end{aligned}$$

By virtue of Assumption 4, each numerator in (51a) and the constraint functions in (51c) are the difference of concave functions, and therefore are not concave in general. In principle, the same approach used for the GEE could be used here to obtain a power allocation algorithm with affordable complexity. This is formalized in the following proposition.

**Proposition 7.** *For any given $\mathbf{p}_j$, denote by $\mathcal{G}_j$ the optimiza-*

[4] We observe that Slater's condition is a suitable constraint qualification for pseudo-concave maximization problems subject to convex constraints [47].

[5] Recall that $\dfrac{\partial \log_2 |\mathbf{A} + x\mathbf{B}|}{\partial x} = \text{tr}\left((\mathbf{A} + x\mathbf{B})^{-1}\mathbf{B}\right)$, for a scalar $x$ and $\mathbf{A}, \mathbf{B}$ being square matrices of proper dimensions [48].

*tion problem*

$$\max_{(t,\mathbf{p})} \min_{k=1,\ldots,K} \frac{q_k^+(\mathbf{p}) - \left[q_k^-(\mathbf{p}_j) + \left(\nabla_{\mathbf{p}} q_k^-|_{\mathbf{p}=\mathbf{p}_j}\right)^T (\mathbf{p}-\mathbf{p}_j)\right]}{\mu_k p_k + \Psi_k} \tag{52a}$$

$$\text{s.t. } 0 \leq p_k \leq P_{\max,k} \ \forall k \tag{52b}$$

$$c_k^+(\mathbf{p}) - \left[c_k^-(\mathbf{p}_j) + \left(\nabla_{\mathbf{p}} c_k^-|_{\mathbf{p}=\mathbf{p}_j}\right)^T (\mathbf{p}-\mathbf{p}_j)\right] \geq 0 \ \forall k \tag{52c}$$

*and call $\mathbf{p}_j^\star$ its optimal solution. If $\mathbf{p}_j = \mathbf{p}_{j-1}^\star \ \forall j \geq 1$, and $\mathbf{p}_0$ is any feasible power vector, then $\{\text{WMEE}(\mathbf{p}_j^\star)\}_j$ is monotonically increasing and converges.*

*Proof:* The proof follows along the same line of reasoning used for Proposition 6. ■

Since, the numerator and denominator of (52a) are concave and convex, respectively, and all constraint functions are concave or affine, it is possible to globally solve (52) by means of fractional programming theory with affordable complexity. Also, it should be mentioned that KKT-optimality is not mentioned in the statement of Proposition 7, for the simple reason that, unlike the GEE, the WMEE is not differentiable, thus implying that (51) admits no KKT conditions. However, it is possible to reformulate (51) through its equivalent epigraph-form representation [38], expressed as:

$$\underset{(t,\mathbf{p})}{\text{maximize}} \quad t \tag{53a}$$

$$\text{s.t.} \quad 0 \leq p_k \leq P_{\max,k} \ \forall k \tag{53b}$$

$$c_k^+(\mathbf{p}) - c_k^-(\mathbf{p}) \geq 0 \ \forall k \tag{53c}$$

$$q_k^+(\mathbf{p}) - q_k^-(\mathbf{p}) - t(\mu_k p_k + \Psi_k) \geq 0 \ \forall k. \tag{53d}$$

Then, it is possible to apply the sequential method to obtain a first-order optimal solution of (53).

**Remark 3.** Similarly as for Section IV, the approaches described in this section can be extended to beamforming problems, too. This would result again in optimization problems with objectives and constraints expressed as the difference of concave functions, with the additional difficulty that the functions would be matrix-variate. However, the approach developed here extends to this case, because it is always possible to linearize the non-concave part of the numerators of the energy efficiencies with respect to the matrix variables.

## VI. Numerical Examples

Among the many possible scenarios under the umbrella of the proposed optimization frameworks, we focus on the two case-studies of a multi-antenna LTE network, the leading standard in 4G systems, and of a massive MIMO network, one of the strongest candidate technologies for 5G networks.

### A. MIMO LTE network with LMMSE detection

Consider the uplink of a multi-cell LTE system in which the same resource block is used by multiple user equipments (UEs). Each UE and each base station (BS) are equipped with $N_T$ and $N_R$ antennas, respectively. Let us denote by $\mathbf{H}_{k,\ell} \in \mathbb{C}^{N_R \times N_T}$ the channel matrix between UE $k$ and BS $\ell$, while $\mathbf{b}_k \in \mathbb{C}^{N_T}$ is the unit-norm beamforming vector, and $s_k$ is the unit-modulus information symbol sent by UE $k$. UE $k$ is associated with BS $a(k)$. In order to perform data detection for user $k$, the signal received at BS $a(k)$ is linearly processed by a filter $\mathbf{c}_k$. Thus, the signal received at BS $a(k)$ is

$$\mathbf{y}_{a(k)} = \mathbf{c}_k^H \left( \sqrt{p_k} \mathbf{H}_{k,a(k)} \mathbf{b}_k s_k + \sum_{i=1, i\neq k}^{K} \sqrt{p_i} \mathbf{H}_{i,a(k)} \mathbf{b}_i s_i + \mathbf{n}_{a(k)} \right) \tag{54}$$

with $\mathbf{n}_{a(k)} \sim \mathcal{CN}(0,\sigma^2)$ modeling the receiver noise. Assuming $\mathbf{c}_k$ is the LMMSE detector, which is known to be the optimal linear receive structure [29], the SINR enjoyed by UE $k$ turns out to be formally equivalent to (8), with $\mathbf{v}_{ki} = \mathbf{H}_{i,a(k)} \mathbf{b}_i$, $\mathbf{u}_k = \mathbf{0}$, and $r = N_R$. Thus, the functions $q_k^+(\mathbf{p})$ and $q_k^-(\mathbf{p})$ are expressed as in (9) and (10).

In our numerical simulations, we considered a multi-cell system which covers a square area of $2\,\text{km} \times 2\,\text{km}$, wherein UEs are randomly placed and equipped with $N_T = 2$ antennas each. The area is served by $L = 3$ BSs placed at coordinates $(0.5; 0.5)\,\text{km}$, $(0.5; -0.5)\,\text{km}$, $(-0.5; 0)\,\text{km}$, with respect to a reference system with the origin at the center of the square, and UEs are associated to the nearest BS. All propagation channels are generated as realizations of uncorrelated Rayleigh fading, using the path-loss model in [49] with power decay factor equal to 3.5. All mobiles have the same maximum feasible power $P_{\max} = -20\,\text{dBW}$ and hardware-dissipated power $\Psi_k = -20\,\text{dBW}$. The receiver noise power is generated as $\sigma^2 = FB\mathcal{N}_0$, wherein $F = 3\,\text{dB}$ is the noise figure, $B = 180\,\text{Hz}$ is the communication bandwidth, and $\mathcal{N}_0 = -174\,\text{dBm/Hz}$ is the thermal noise power spectral density.

Fig. 1 shows the energy-efficient Pareto region of the system for $K = 2$ UEs. 200 sample points at the Pareto boundary are obtained by solving the WMEE maximization problem for 200 different choices of the weights, each corresponding to finding the outmost point in a certain search direction. The maximization problem was solved using the proposed monotonic fractional programming framework. As a comparison, we show the non-uniform grid of operating points that are achieved by a grid search over 40000 equally-spaced feasible transmit power points. We note that the monotonic fractional programming framework is able to characterize the complete region, while the 40000 points from the grid search fail to explore all parts of the region. Two particular operating points are shown as reference. The point where both UEs use full power is in the interior of the region and thus inefficient, for any choice of EE metric. The maximum GEE points found by sequential fractional programming and monotonic fractional programming coincide, and are also in the interior of the region. This is an interesting result, as it shows that GEE maximization does not necessarily yield a point on the boundary of the energy-efficient Pareto region, thus showing how the GEE metric might fail to capture the efficiency of the individual links. A similar scenario is illustrated in Fig. 2, which considers the case with $K = 3$ UEs. In this case, the Pareto region is illustrated using a number of lines that lie on the boundary. All star-marked points on each line were computed using the proposed monotonic fractional programming framework. It is interesting to observe that the obtained region does not define a convex set, which implies that in general only WMEE maximization can guarantee to characterize the energy-efficient Pareto region of wireless networks.

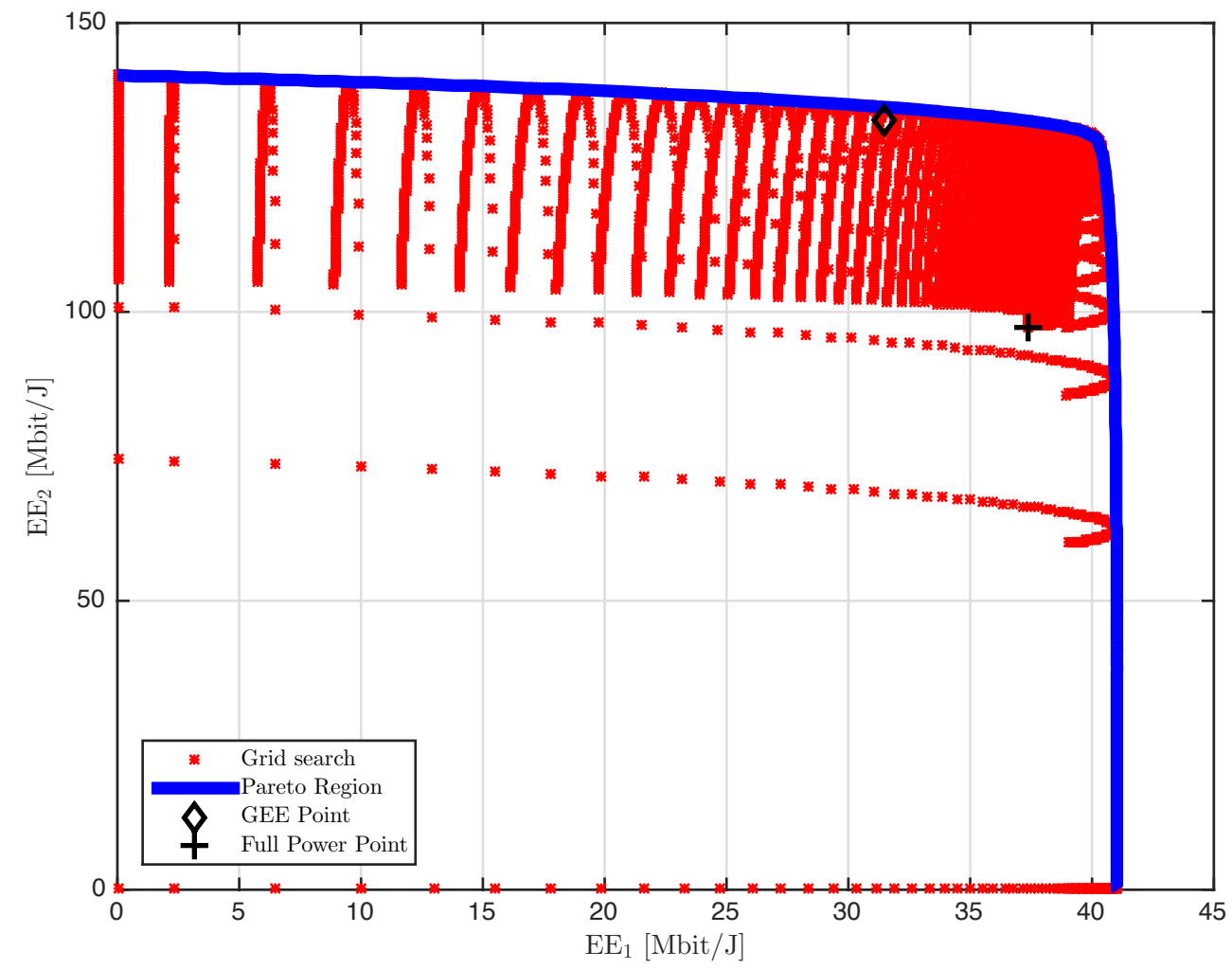

Fig. 1. Energy-efficient Pareto region for $K = 2$, generated by: 1) Monotonic Fractional Programming; 2) Grid Search. The points corresponding to the maximum GEE and to Full Power Allocation are also reported, which are both inside the Pareto region.

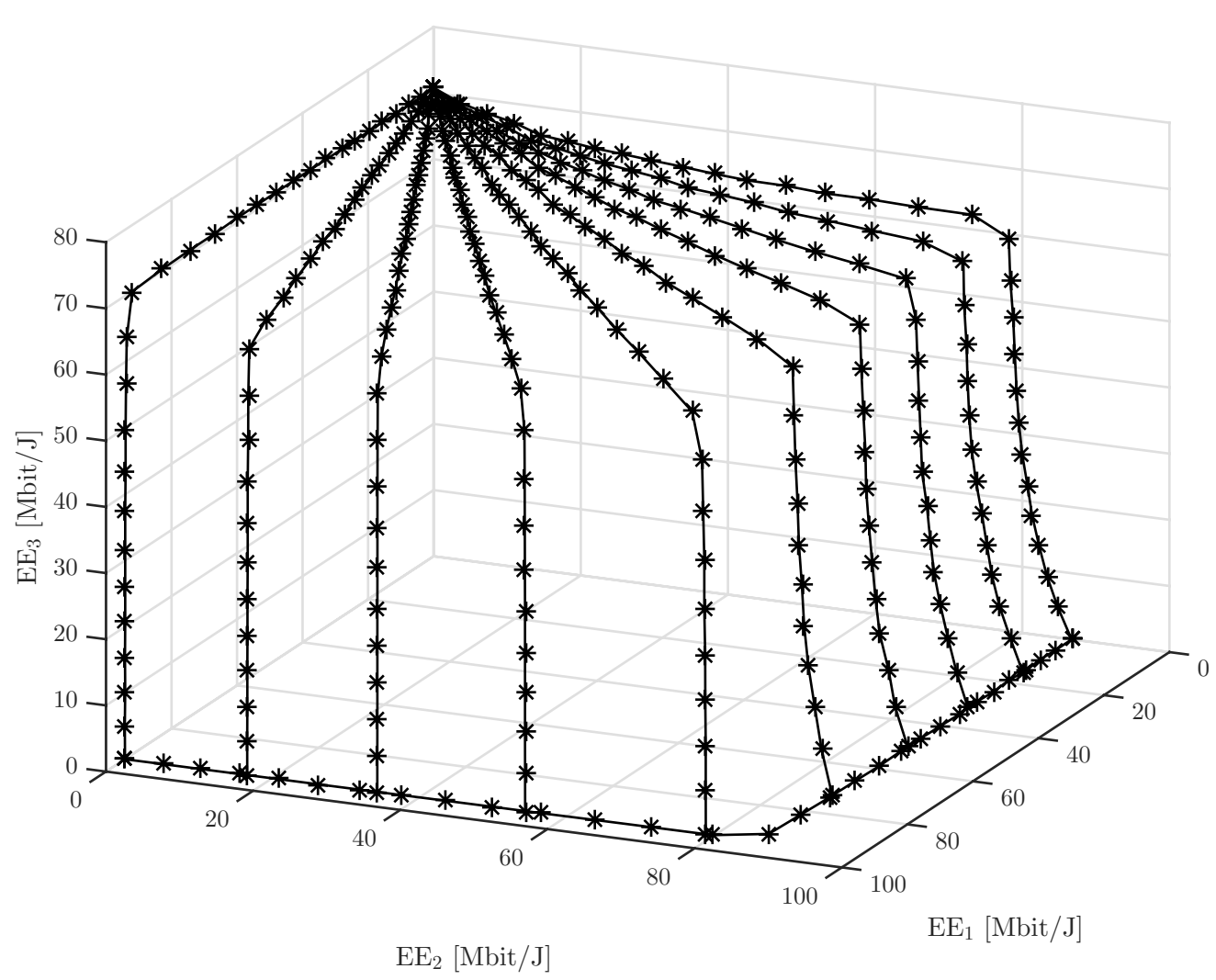

Fig. 2. Energy-efficient Pareto region for $K = 3$, generated by Monotonic Fractional Programming. The region can be seen to define a non-convex set.

### *B. Massive multiple-input multiple-output (MIMO) networks*

Consider the uplink[6] of a massive MIMO network composed of $K$ single-antenna UEs and $L$ BSs (which can represent either macro cells or small cells), each equipped with $N_R$ antennas. Each UE $k$ is associated to a specific serving BS $a(k)$, while interfering with all other UEs. Let us define by $\mathbf{h}_{k,m} \in \mathbb{C}^{N_R}$ the channel vector from UE $k$ to BS $m$. Denoting by $\mathbf{c}_k \in \mathbb{C}^{N_R}$ the linear detector that BS $a(k)$ applies to its received signal to detect the signal from UE $k$, a lower bound[7] on the uplink SINR of UE $k$ takes the following form [50]:

$$\gamma_k = \frac{p_k \left| \mathrm{E}\{\mathbf{c}_k^H \mathbf{h}_{k,a(k)}\} \right|^2}{p_k \mathrm{var}\{\mathbf{c}_k^H \mathbf{h}_{k,a(k)}\} + \mathbf{i}_k} \tag{55}$$

with $\mathbf{i}_k = \sigma^2 \mathrm{E}\{\|\mathbf{c}_k\|^2\} + \sum_{i=1,i\neq k}^{K} p_i \mathrm{E}\{|\mathbf{c}_k^H \mathbf{h}_{i,a(k)}|^2\}$. Assume that a maximum ratio combining (MRC) detector is employed. This amounts to setting $\mathbf{c}_k = \hat{\mathbf{h}}_{k,a(k)}$ where $\hat{\mathbf{h}}_{k,a(k)}$ denotes the estimate of $\mathbf{h}_{k,a(k)}$ given by $\mathbf{h}_{k,a(k)} = \hat{\mathbf{h}}_{k,a(k)} + \tilde{\mathbf{h}}_{k,a(k)}$, with $\tilde{\mathbf{h}}_{k,a(k)}$ being the estimation error statistically independent of $\hat{\mathbf{h}}_{k,a(k)}$. We consider Rayleigh fading channels $\mathbf{h}_{k,m} \sim \mathcal{CN}(0, d_{k,m}\mathbf{I}_{N_R})$ where the variance $d_{k,m}$ accounts for the large-scale channel fading and attenuation from UE $k$ to BS $m$. If a minimum mean square error (MMSE)-based channel estimation scheme is used at the BS (with full pilot reuse across cells) [50], then we have that $\hat{\mathbf{h}}_{k,a(k)} \sim \mathcal{CN}(\mathbf{0}, \rho_{k,a(k)}\mathbf{I}_{N_R})$ and $\tilde{\mathbf{h}}_{k,a(k)} \sim \mathcal{CN}(\mathbf{0}, (d_{k,a(k)} - \rho_{k,a(k)})\mathbf{I}_{N_R})$ where $\rho_{k,a(k)} = \frac{d_{k,a(k)}}{\tau + \sum_m d_{k,m}}$, with $\tau$ being a given parameter that depends on the total pilot transmit power over the pilot sequence. Under the above assumptions, we have that $\gamma_k$ takes the form in (5), with $\alpha_k = \rho_{k,a(k)}$, $\phi_k = d_{k,a(k)} + \sum_{m\neq a(k)} \rho_{k,m}^2/\rho_{k,a(k)}$, and $\beta_{i,k} = d_{i,a(k)}\rho_{i,a(k)}/\rho_{k,a(k)}$, which is formally equivalent to (5). Thus, the functions $q_k^+(\mathbf{p})$ and $q_k^-(\mathbf{p})$ are expressed as in (6) and (7) for all $k = 1, \ldots, K$.

Similar results can be obtained when the system is affected by hardware impairments [51], [52]; for example, unavoidable clock drifts in local oscillators, finite-precision digital-to-analog converters, amplifier non-linearities, non-ideal analog filters, and so forth. For the sake of simplicity, let us assume that the hardware impairments are only present at the UEs.[8] Following [51], hardware impairments result in a reduction of the uplink signals by a factor $1 - \epsilon^2$, with $\epsilon$ being the error vector magnitude, and in an additive Gaussian distortion noise which carries the removed useful power. In these circumstances, a lower bound of the SINR is given by

$$\gamma_k = \frac{p_k(1-\epsilon^2) \left| \mathrm{E}\{\mathbf{c}_k^H \mathbf{h}_{k,a(k)}\} \right|^2}{p_k(1-\epsilon^2)\mathrm{var}\{\mathbf{c}_k^H \mathbf{h}_{k,a(k)}\} + p_k\epsilon^2 \mathrm{E}\{|\mathbf{c}_k^H \mathbf{h}_{k,a(k)}|^2\} + \mathbf{i}_k} \tag{56}$$

Plugging $\mathbf{c}_k = \hat{\mathbf{h}}_{k,a(k)}$ into the above equation and taking into account that in the presence of hardware impairments $\hat{\mathbf{h}}_{k,a(k)} \sim \mathcal{CN}(\mathbf{0}, \sqrt{1-\epsilon^2}\rho_{k,a(k)}\mathbf{I}_{N_R})$ and $\tilde{\mathbf{h}}_{k,a(k)} \sim \mathcal{CN}(\mathbf{0}, (d_{k,a(k)} - \sqrt{1-\epsilon^2}\rho_{k,a(k)})\mathbf{I}_{N_R})$, the SINR is easily found to be again in the same form of (5).

In our numerical simulations, we consider a cellular setup wherein 3 small-cell receivers equipped with 20 antennas each are deployed in a similar way as in Section VI-A. In addition, one macro-BS with 50 receive antennas is placed at the center of the area to serve. Instead, all mobiles have a single antenna. We consider the presence of hardware impairments with $\epsilon = 10^{-1}$, and of channel estimation errors with $\tau = 0.3$.

[6]Similar results can be obtained for the downlink, but we selected the uplink to reduce the amount of notation.

[7]This refers to the standard worst-case lower bound on the mutual information where the uncorrelated interference is treated as Gaussian noise, which is information-theoretic optimal for small interference powers.

[8]The impact of BS hardware impairments vanishes as $N_R$ increases [52], thus UE hardware impairments are expected to dominate in massive MIMO networks.

Our numerical experiments confirm that the sequential fractional programming framework performs as the monotonic fractional programming framework and thus achieves the global GEE maximum, and not a suboptimal solution. Recall that the latter framework is guaranteed to find the GEE maximum, but with a computational complexity that grows exponentially in the number of UEs. If the algorithm is initiated at $\lambda_0 = 0$, corresponding to zero GEE, it finds the sum-rate maximizing solution in the first iteration of Dinkelbach's algorithm. In our experiments, only one or two further iterations are required to converge to the global GEE optimum, which is line with the super-linear convergence rate of Dinkelbach's algorithm. At convergence, the difference between the current GEE value $\lambda_j$ and the next obtained GEE value $\lambda_{j+1}$ is negligible. Fig. 3 shows the behavior when $\lambda_j$ equals the GEE value obtained by the sequential fractional programming framework with $K = 2$. Fig. 4 shows the corresponding behavior for $K = 3$. In both cases, the lower and upper bounds in the BRB algorithm that solves the monotonic subproblem converge to the same GEE value ($\lambda_{j+1} = \lambda_j$), which validates that the algorithm has already converged to the global maximum. The bounding behavior is very typical for the BRB algorithm [27]; a relatively small difference between the lower and upper bounds is obtained quickly, while many more iterations in the BRB algorithm are required to push the difference down to zero. Notice that 40 iterations are sufficient for $K = 2$, while the number of iterations for $K = 3$ is of the order of $10^4$, which shows the exponential complexity scaling with the number of UEs.

Finally, Figs. 5 and 6 compare the GEE achieved by the monotonic fractional programming and the sequential fractional programming frameworks, versus $P_{\max}$, for $K = 2$ and $K = 3$, respectively. The GEE obtained when the sequential framework is used to maximize the sum-rate, and when all users transmit with full power (i.e. $p_k = P_{\max}$ for all $k$) are shown, too. The results again confirm the optimality of the sequential approach, which performs as the monotonic approach. In addition, both schemes saturate at large $P_{\max}$, because once $P_{\max}$ is large enough to allow attaining the GEE global maximum, the excess transmit power is no longer used. Using full power at all UEs is globally optimal up to $P_{\max} = -34\,\text{dBW}$, whereas using all the excess power at larger $P_{\max}$ will degrade the GEE. The GEE obtained by sum-rate maximization is globally optimal up to $P_{\max} = -28\,\text{dBW}$, while it decreases for larger $P_{\max}$. This shows that at least one UE should transmit at $-28\,\text{dBW}$ at the globally optimal GEE point, but not all of them.

Next, we provide numerical results when a larger number of mobiles is to be optimized. In particular, we assume $K = 8$ mobiles are present in the area to cover, and Fig. 7 compares the GEE achieved by the monotonic fractional programming and the sequential fractional programming frameworks, versus $P_{\max}$. Again, the GEE obtained by sum-rate maximization and by full power transmission are reported for comparison purposes. The results again indicate that the sequential fractional programming framework performs as the globally optimal monotonic fractional programming framework. Moreover, similar remarks as for previous illustrations can be made.

Finally, considering the same scenario as in Fig. 7, Tab. VI-B shows the number of iterations to reach convergence for the monotonic fractional programming approach and for the sequential fractional programming approach versus $P_{max}$. For the former algorithm, the reported value is just the number of iterations required by Dinkelbach's algorithm to converge, whereas for the latter approach the reported value is the number of outer iterations, i.e. the number of approximate fractional problems to be solved by fractional programming. In both cases the convergence was declared when the relative squared error between the GEE values between two successive iterations was not larger than $10^{-4}$. Full power allocation was used as initialization point. The results show that the monotonic-based algorithm converges in a very limited number of iterations. This is expected since the monotonic approach allows for the optimal implementation of Dinkelbach's method, which is known to have a super-linear convergence. Instead, a slightly larger, but still limited, number of iterations is required for the sequential method. However, it is important to remember that solving each approximate problem in the sequential method can be accomplished with limited complexity, whereas the globally optimal method requires to solve a monotonic problem in each Dinkelbach iteration. Thus, the results show how the sequential method indeed lends itself to being implemented in practical systems, while the monotonic approach is useful for benchmarking purposes. Finally, it is seen that the number of iterations increases with $P_{max}$, which is expected since a larger $P_{max}$ means a larger feasible set over which to optimize.

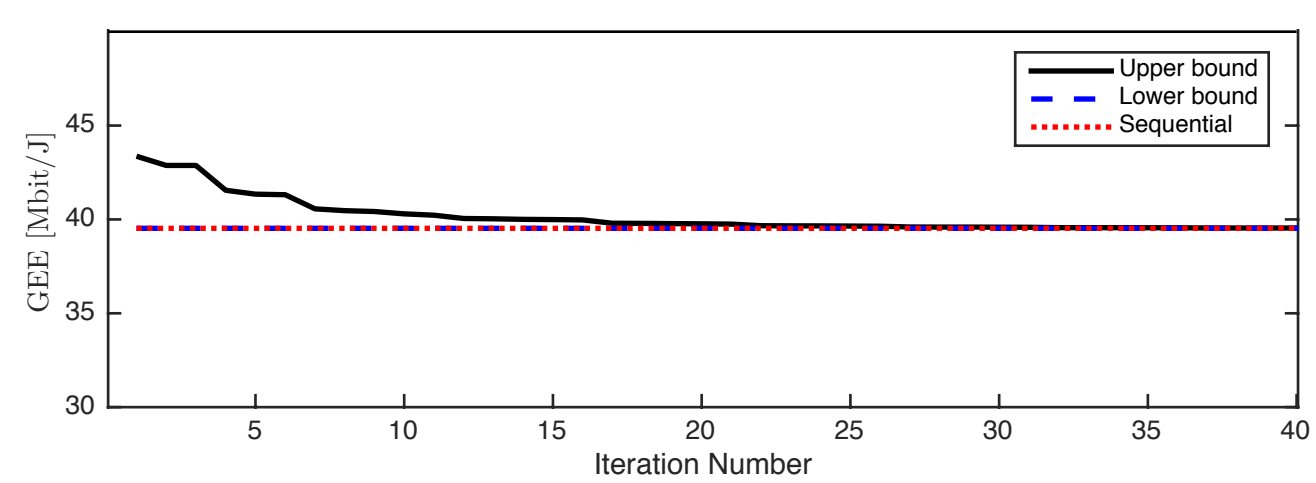

Fig. 3. Convergence behavior of the BRB algorithm, with $K = 2$, in the last iteration of Dinkelbach's algorithm. The behavior shows the convergence to the global GEE optimal solution and illustrates the optimality of the sequential fractional programming framework.

## VII. Conclusion

This work has developed two optimization frameworks to characterize the ultimate energy-efficient performance of wireless networks. The former merges the tool of monotonic optimization with the theory of fractional programming, and is guaranteed to achieve global optimality, with a complexity that, although still exponentially increasing with the number of links, is significantly lower than that of standard global optimization methods. The latter combines the tool of sequential optimization with fractional programming theory and enjoys first-order optimality with affordable complexity. Numeri-

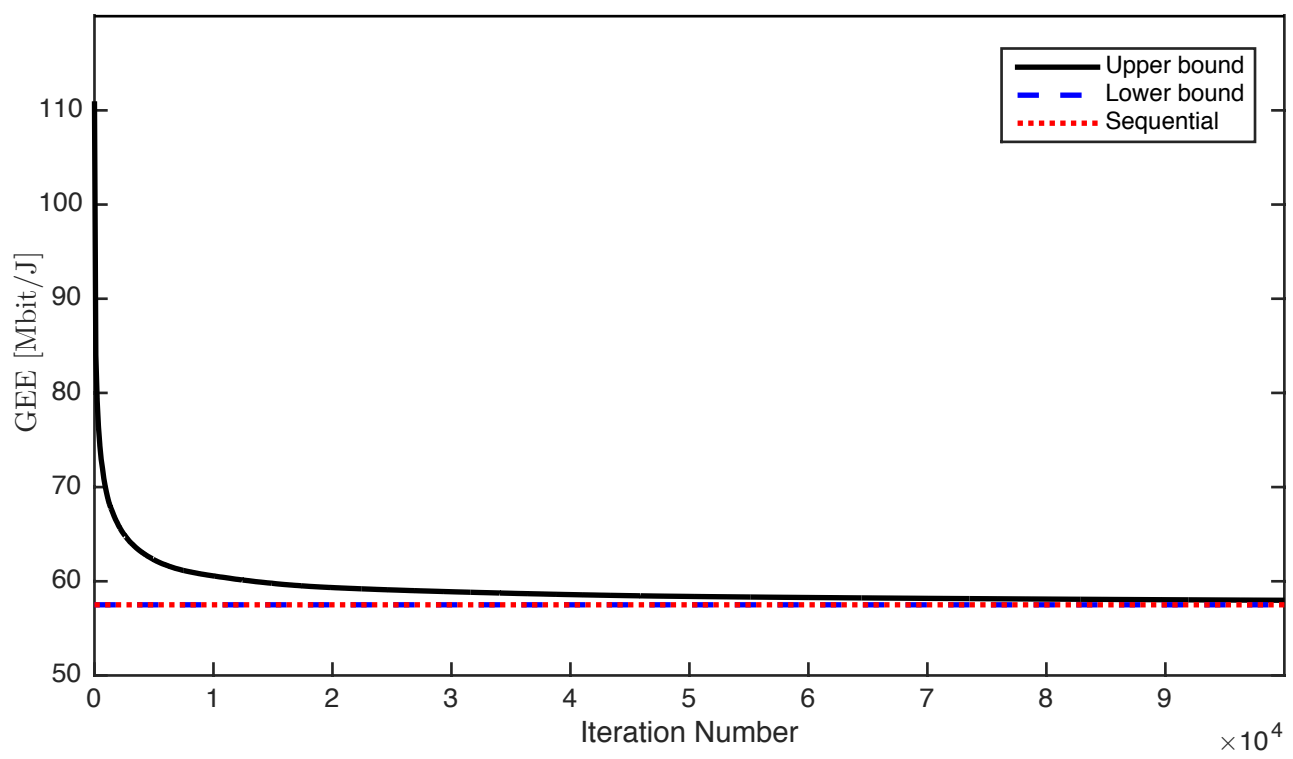

Fig. 4. Convergence behavior of the BRB algorithm, with $K = 3$, in the last iteration of Dinkelbach's algorithm. The behavior shows the convergence to the global GEE optimal solution and illustrates the optimality of the sequential fractional programming framework.

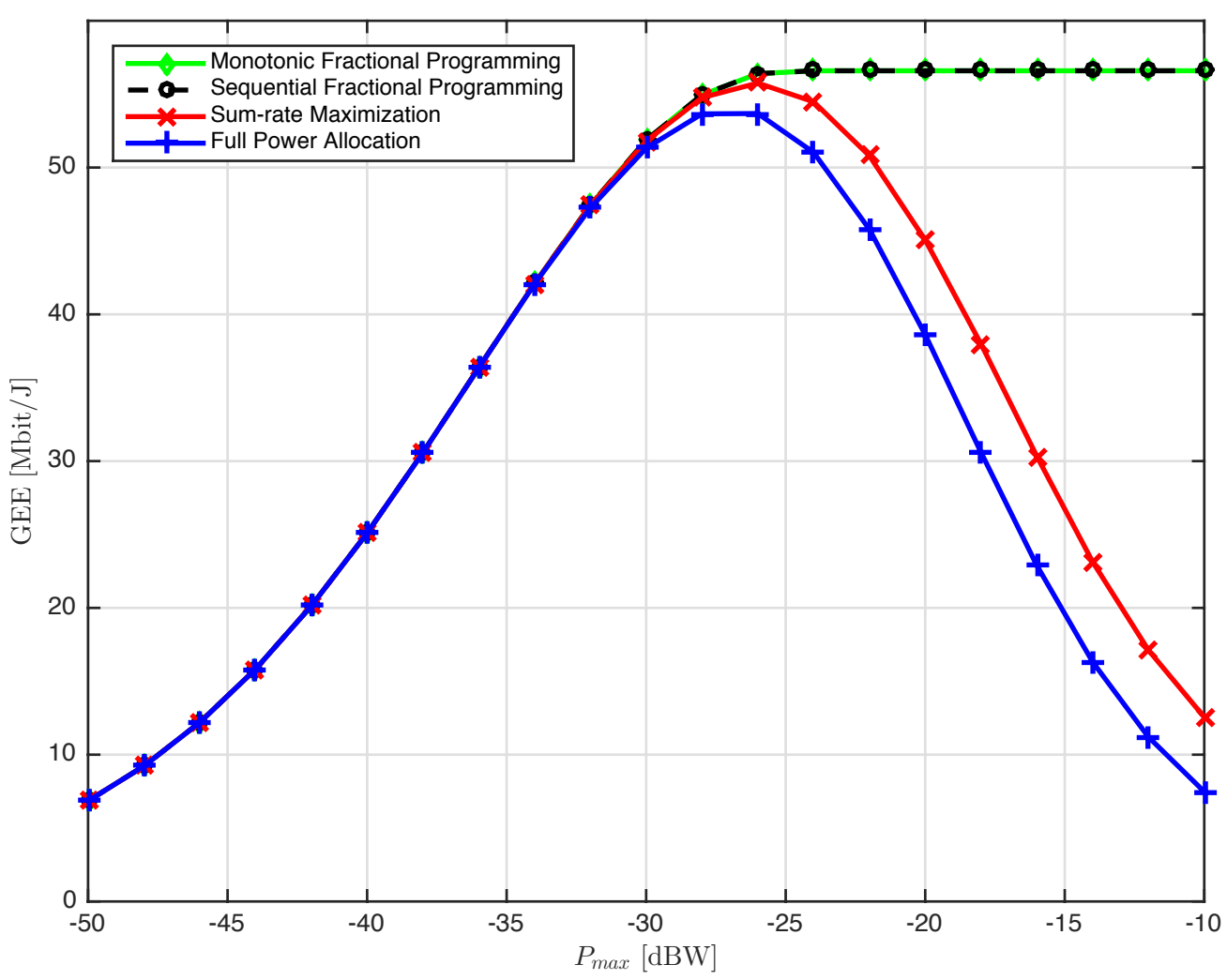

Fig. 5. Achieved GEE versus $P_{\max}$ for $K = 2$, using: 1) Monotonic fractional programming; 2) Sequential Fractional Programming, 3) Sum-rate maximization by sequential programming; 4) Full Power Allocation.

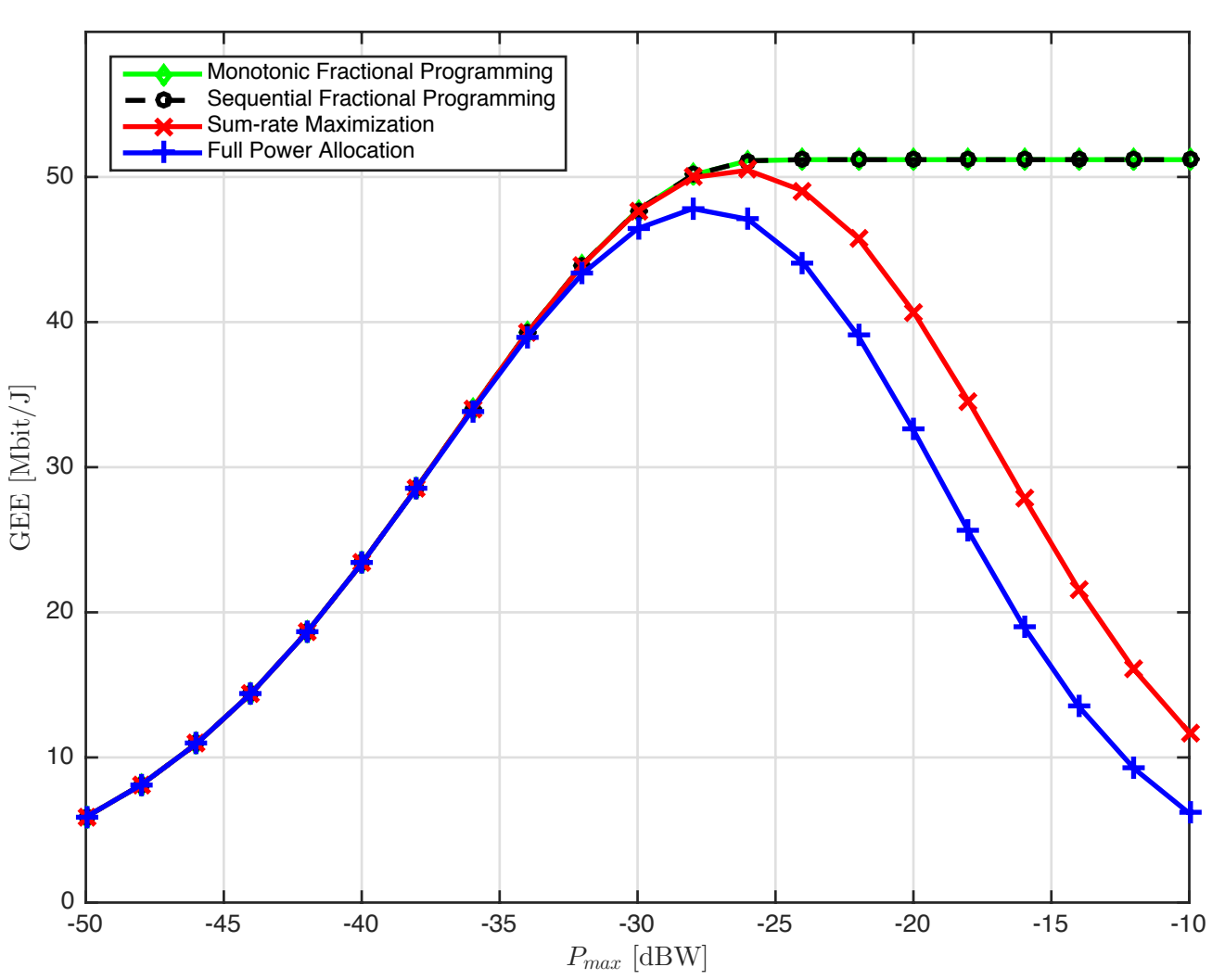

Fig. 6. Achieved GEE versus $P_{\max}$ for $K = 3$, using: 1) Monotonic fractional programming; 2) Sequential Fractional Programming, 3) Sum-rate maximization by sequential programming; 4) Full Power Allocation.

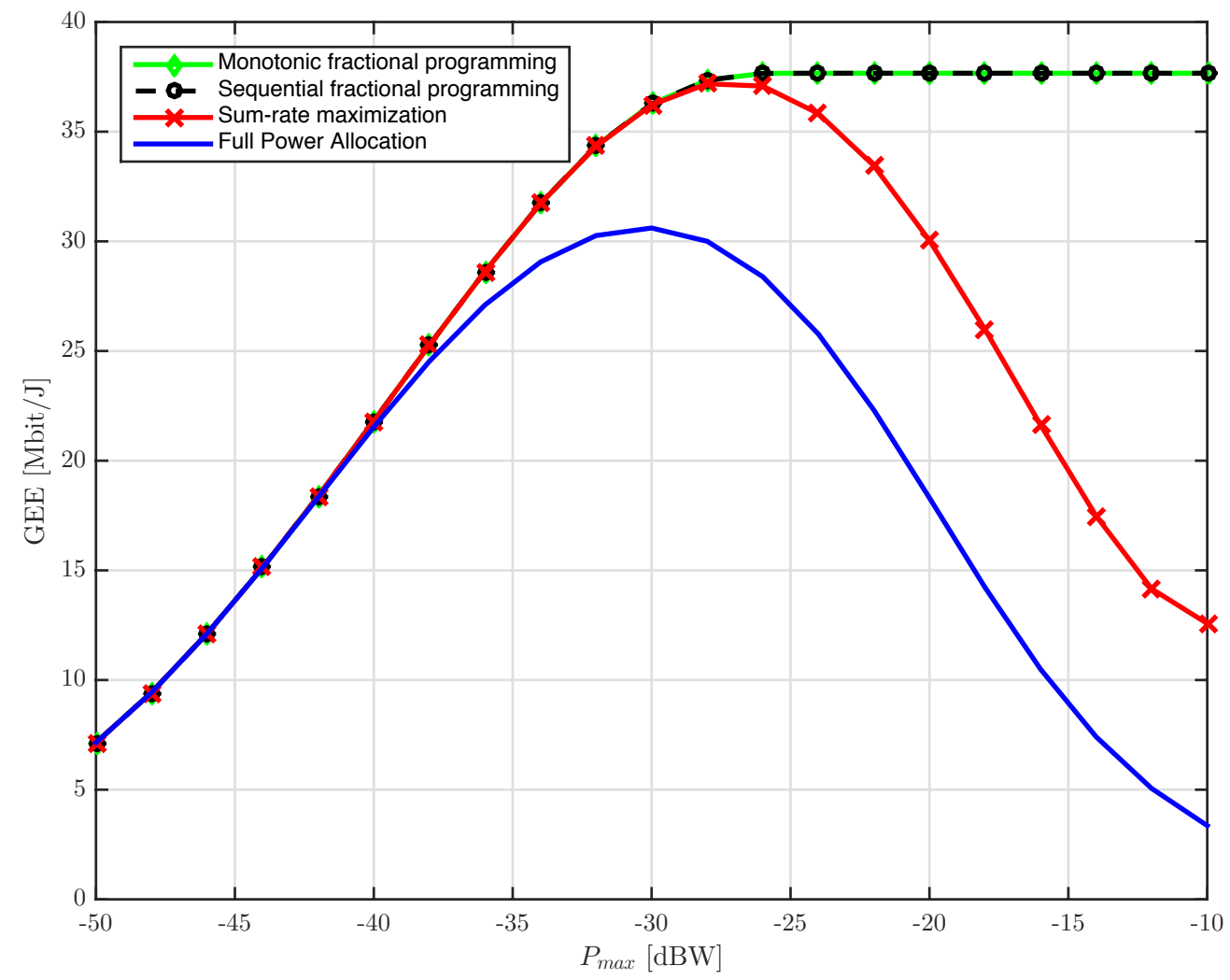

Fig. 7. Achieved GEE versus $P_{\max}$ for $K = 8$, using: 1) Monotonic fractional programming; 2) Sequential Fractional Programming, 3) Sum-rate maximization by sequential programming; 4) Full Power Allocation.

| $P_{max}$ | Monotonic | Sequential |
|---|---|---|
| $-50$ dBW | 1 | 1 |
| $-45$ dBW | 1 | 1 |
| $-40$ dBW | 1.57 | 2.17 |
| $-35$ dBW | 2.23 | 3.01 |
| $-30$ dBW | 2.47 | 4.17 |
| $-25$ dBW | 2.81 | 5.90 |
| $-20$ dBW | 3.01 | 7.77 |
| $-15$ dBW | 3.15 | 8.56 |
| $-10$ dBW | 3.28 | 8.64 |

TABLE I

AVERAGE NUMBER OF ITERATIONS VERSUS $P_{max}$ FOR $K = 8$, USING: 1) MONOTONIC FRACTIONAL PROGRAMMING (I.E. OPTIMAL IMPLEMENTATION OF DINKELBACH'S ALGORITHM); 2) SEQUENTIAL FRACTIONAL PROGRAMMING.

cal results also show that sequential fractional programming achieves global optimality in several practical scenarios.

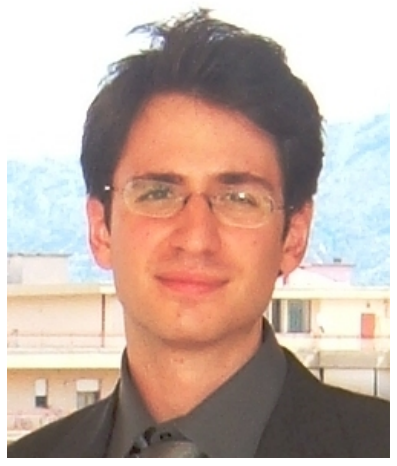

**Alessio Zappone** (S'08 – M'11 – SM'16) received his Ph.D. from the University of Cassino and Southern Lazio, Cassino, Italy. Afterwards, he worked with Consorzio Nazionale Interuniversitario per le Telecomunicazioni (CNIT) in the framework of the FP7 EU-funded project TREND. From 2012 to 2016, Alessio has been with the Department of Communication Theory of the Technische Universität Dresden, Dresden, Germany, serving as the principal investigator of the project CEMRIN, funded by the German research foundation (DFG) and carried out at the Department of Communication Theory of the Technische Universität Dresden, Dresden, Germany. Since 2016 he is adjunct professor with the University of Cassino and Southern Lazio.

His research interests lie in the area of communication theory and signal processing, with main focus on optimization techniques for resource allocation and energy efficiency. He held several research appointments at TU Dresden, Politecnico di Torino, Suplec - Alcatel-Lucent Chair on Flexible Radio, and University of Naples Federico II. He was the recipient of a Newcom# mobility grant in 2014. Alessio serves as associate editor for the IEEE Signal Processing Letters and has served as associate editor for the IEEE Journal on Selected Areas on Communications (Special Issue on Energy-Efficient Techniques for 5G Wireless Communication Systems). He is the first author of the textbook *Energy efficiency in wireless networks via fractional programming theory* (Foundations and Trends in Communications and Information Theory, 2015).

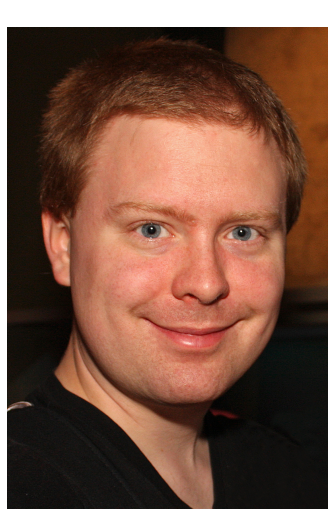

**Emil Björnson** (S'07, M'12) received his M.S. degree in Engineering Mathematics from Lund University, Sweden, in 2007. He received his Ph.D. degree in Telecommunications from the KTH Royal Institute of Technology, Stockholm, Sweden, in 2011. From 2012 to July 2014, he was a joint postdoc at Suplec, Gif-sur-Yvette, France, and at KTH Royal Institute of Technology. He is currently an Assistant Professor at the Department of Electrical Engineering (ISY) at Linköping University, Sweden.

His research interests include multi-antenna cellular communications, radio resource allocation, energy efficiency, massive MIMO, and network topology design. He is on the editorial board of *IEEE Transactions on Green Communications and Networking* since 2016. He is the first author of the textbook *Optimal Resource Allocation in Coordinated Multi-Cell System* (Foundations and Trends in Communications and Information Theory, 2013). He is also dedicated to reproducible research and has made a large amount of simulation code publicly available. Dr. Björnson received the 2014 Outstanding Young Researcher Award from IEEE ComSoc EMEA, the 2015 Ingvar Carlsson Award, and the 2016 Best PhD Award. He has received five best paper awards for novel research on optimization and design of multi-cell multi-antenna communications: ICC 2015, WCNC 2014, SAM 2014, CAMSAP 2011, and WCSP 2009.

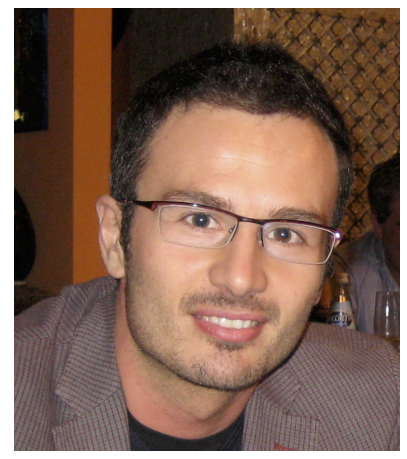

**Luca Sanguinetti** (SM'15) received the Laurea Telecommunications Engineer degree (cum laude) and the Ph.D. degree in information engineering from the University of Pisa, Italy, in 2002 and 2005, respectively. Since 2005 he has been with the Dipartimento di Ingegneria dell'Informazione of the University of Pisa. In 2004, he was a visiting Ph.D. student at the German Aerospace Center (DLR), Oberpfaffenhofen, Germany. During the period June 2007 - June 2008, he was a postdoctoral associate in the Dept. Electrical Engineering at Princeton. During the period June 2010 - Sept. 2010, he was selected for a research assistantship at the Technische Universität Munchen. From July 2013 to July 2015, he was with the Alcatel-Lucent Chair on Flexible Radio, Supélec, Gif-sur-Yvette, France. He is an Assistant Professor at the Dipartimento di Ingegneria dell'Informazione of the University of Pisa. L. Sanguinetti is currently serving as an Associate Editor for the IEEE Transactions on Wireless communications, IEEE Signal Processing Letters. He served as Lead Guest Editor of IEEE Journal on Selected Areas of Communications Special Issue on "Game Theory for Networks" and as an Associate Editor for IEEE Journal on Selected Areas of Communications (series on Green Communications and Networking). His expertise and general interests span the areas of communications and signal processing, game theory and random matrix theory for wireless communications. He was the co-recipient of two best paper awards: *IEEE Wireless Commun. and Networking Conference (WCNC) 2013* and *IEEE Wireless Commun. and Networking Conference (WCNC) 2014*. He was also the recipient of the FP7 Marie Curie IEF 2013 "Dense deployments for green cellular networks".

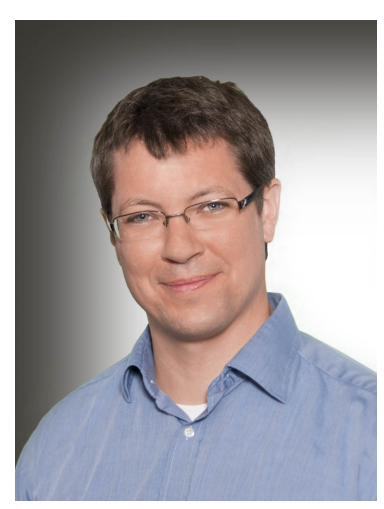

**Eduard Jorswieck** (S'01 – M'03 – SM'08) received the Diplom-Ingenieur (M.S.) degree and Doktor-Ingenieur (Ph.D.) degree, both in electrical engineering and computer science, from the Technische Universität Berlin, Germany, in 2000 and 2004, respectively. He was with the Broadband Mobile Communication Networks Department, Fraunhofer Institute for Telecommunications, Heinrich-Hertz-Institut, Berlin, from 2000 to 2008. From 2005 to 2008, he was a Lecturer with the Technische Universität Berlin. From 2006 to 2008, he was with the Department of Signals, Sensors and Systems, Royal Institute of Technology, as a Post-Doctoral Researcher and an Assistant Professor. Since 2008, he has been the Head of the Chair of Communications Theory and a Full Professor with the Technische Universität Dresden, Germany. He is principal investigator in the excellence cluster center for Advancing Electronics Dresden (cfAED) and founding member of the 5G lab Germany (5Glab.de).

His main research interests are in the area of signal processing for communications and networks, applied information theory, and communications theory. He has authored over 80 journal papers, 8 book chapters, some 225 conference papers and 3 monographs on these research topics. Eduard was a co-recipient of the IEEE Signal Processing Society Best Paper Award in 2006 and co-authored papers that won the Best Paper or Best Student Paper Awards at IEEE WPMC 2002, Chinacom 2010, IEEE CAMSAP 2011, IEEE SPAWC 2012, and IEEE WCSP 2012.

Dr. Jorswieck was a member of the IEEE SPCOM Technical Committee (2008 - 2013), and has been a member of the IEEE SAM Technical Committee since 2015. Since 2011, he has been an Associate Editor of the IEEE TRANSACTIONS ON SIGNAL PROCESSING. Since 2008, continuing until 2011, he has served as an Associate Editor of the IEEE SIGNAL PROCESSING LETTERS, and until 2013, as a Senior Associate Editor. Since 2013, he has served as an Editor of the IEEE TRANSACTIONS ON WIRELESS COMMUNICATIONS.